\documentclass[10pt,leqno]{article} 

\usepackage{amssymb} 
\usepackage{amsmath} 

\title{\sc{Automorphisms and forms of simple infinite-dimensional
linearly compact Lie superalgebras}}

\author{{\sc Nicoletta Cantarini}\thanks{Partially supported by
local research grants of the University of Padova}\\
\normalsize{\it Dipartimento di Matematica Pura ed Applicata,}\\ 
\normalsize{\it Universit\`a di Padova, Via Belzoni, 7}\\
\normalsize{\it Padova, 35131, Italy}\\
\normalsize{\it cantarin@math.unipd.it}\vspace{0.3cm}
\and\setcounter{footnote}{6}
{\sc Victor G.\ Kac}\thanks{Partially supported by NSF grant DMS0501395}\\
\normalsize{\it Department of Mathematics,}\\
\normalsize{\it Massachusetts Institute of Technology, 77 Massachusetts Ave}\\
\normalsize{\it Cambridge, Massachusetts 02139, USA}\\
\normalsize{\it kac@math.mit.edu}}

\newtheorem{theorem}{Theorem}[section] 
\newtheorem{lemma}[theorem]{Lemma} 
\newtheorem{corollary}[theorem]{Corollary} 
\newtheorem{proposition}[theorem]{Proposition} 
\newtheorem{definition}[theorem]{Definition} 
\newtheorem{remark}[theorem]{Remark}
\newtheorem{example}[theorem]{Example}
 
\def\Z{\mathbb{Z}} 

\def\F{\mathbb{F}}
\def\g{\mathfrak{g}}

\def\h{\mathfrak{h}} 
\def\N{\mathbb{N}}
\def\C{\mathbb{C}}

\date{} 
\begin{document} 

\markboth{N.\ Cantarini, V.\ G.\ Kac}
{Automorphisms and forms of simple infinite-dimensional
linearly compact Lie superalgebras}

\maketitle 
\begin{center}
{\em To Dmitri V.\ Alekseevski on his 65th birthday}
\end{center}
\date{}
\medskip 
\begin{abstract} We describe the group of continous automorphisms 
of all simple
infinite-dimensional linearly compact Lie superalgebras and use it
in order to classify $\F$-forms of these superalgebras over any
field $\F$ of characteristic zero.

\medskip

\noindent
{\it Keywords:} Linearly compact Lie superalgebra, $\F$-form, Galois cohomology.
\end{abstract}
\section*{Introduction} 
In our paper \cite{CantaK} we classified all maximal open subalgebras
of all simple infinite-dimensional linearly compact Lie superalgebras
$S$ over an algebraically closed field $\bar{\F}$ of characteristic zero,
up to conjugation by the group $G$ of inner automorphisms of the Lie
superalgebra $Der S$ of continuous derivations of $S$. An immediate corollary
of this result is Theorem 11.1 of \cite{CantaK}, which describes,
up to conjugation by $G$, all maximal open subalgebras of $S$, which are
invariant with respect to all inner automorphisms of $S$.
Using this result and an explicit description of $Der S$ (see
\cite[Proposition 6.1]{K} and its corrected version \cite[Proposition
1.8]{CantaK}), we derive the classification of all maximal among
the open subalgebras of $S$, which are $Aut S$-invariant, where
$Aut S$ is the group of all continuous automorphisms of $S$
(Theorem \ref{autSinv}).
Such a subalgebra $S_0$ always exists, and in most of the cases
it is unique (also, in most of the cases it is a subalgebra
of minimal codimension).
Picking a subspace $S_{-1}$ of $S$, which is minimal among $Aut S$-invariant
subspaces, properly containing $S_0$, we can construct the Weisfeiler
filtration (see e.g.\ \cite{CantaK} or \cite{K}). Then it is easy to
see that
\begin{equation}
Aut S={\cal U}\rtimes Autgr S,
\label{0.1}
\end{equation}
where $\cal{U}$ is a normal prounipotent subgroup consisting of
automorphisms of $S$ inducing an identity automorphism of $Gr S$, and
$Autgr S$ is a subgroup of a (finite-dimensional) algebraic group
of automorphisms of $Gr S$, preserving the grading.

We list all the groups $Autgr S$, along with their (faithful) action
on $Gr_{-1}S$, in Table 1. This leads to the following
description of the group $Aut S$:
\begin{equation}
Aut S=Inaut S\rtimes A,
\label{0.2}
\end{equation}
where $Inaut S$ is the subgroup of all inner automorphisms of $S$ and $A$ 
is a closed subgroup of $Autgr S$, listed in Corollary \ref{outer}.

Let $\F$ be a subfield of $\bar{\F}$, whose algebraic closure is
$\bar{\F}$, and fix an $\F$-form $S^{\F}$ of $S$, i.e., a Lie superalgebra
over $\F$, such that $S^\F\otimes_\F\bar{F}\cong S$. Then all
$\F$-forms of $S$, up to isomorphism, are in a bijective
correspondence with $H^1(Gal, Aut S)$, where $Gal=Gal(\bar{\F}/\F)$
(see e.g.\ \cite{R}). Since the first Galois cohomology of a prounipotent
algebraic group is trivial (see e.g.\ \cite{R}),
we conclude, using the cohomology long exact sequence,
that
\begin{equation}
H^1(Gal, Aut S)\hookrightarrow H^1(Gal, Autgr S).
\label{0.3}
\end{equation}

The infinite-dimensional linearly compact simple Lie superalgebras
have been classified in \cite{K}. The list consists of ten series
$(m\geq 1)$: $W(m,n)$, $S(m,n)$ $((m,n)\neq (1,1))$, $H(m,n)$ ($m$ even), 
$K(m,n)$ ($m$ odd), $HO(m,m)$ $(m\geq 2)$, $SHO(m,m)$ $(m\geq 3)$,
$KO(m,m+1)$, $SKO(m,m+1;\beta)$ $(m\geq 2)$, $SHO^\sim(m,m)$ ($m$ even),
$SKO^\sim(m,m+1)$ ($m\geq 3$ odd), and five exceptional Lie
superalgebras: $E(1,6)$, $E(3,6)$, $E(3,8)$, $E(4,4)$, $E(5,10)$.
Since the following isomorphisms hold (see \cite{CantaK}, \cite{K}):
$W(1,1)\cong K(1,2)\cong KO(1,2)$, $S(2,1)\cong HO(2,2)\cong SKO(2,3;0)$,
$SHO^\sim(2,2)\cong H(2,1)$,
when dealing with $W(m,n)$, $KO(n,n+1)$,
$HO(n,n)$, $SKO(2,3;\beta)$ and $SHO^\sim(n,n)$,
 we will assume that $(m,n)\neq (1,1)$, $n\geq 2$, $n\geq 3$, $\beta\neq 0$,
and $n>2$, respectively. We will use the construction of all these
superalgebras as given in \cite{CantaK} (see also \cite{CCK}, \cite{CK},
\cite{K}, \cite{Gafa}, \cite{S}).

Since the first Galois cohomology with coefficients in the groups
$GL_n(\bar{F})$ and $Sp_n(\bar{\F})$ is trivial (see, e.g.,
\cite{Serre1}, \cite{Serre2}), we conclude
from (\ref{0.3}) and Table 1, that $H^1(Gal, Aut S)$ is trivial in
all cases except for four:
$S=H(m,n)$, $K(m,n)$, $S(1,2)$, and $E(1,6)$. Thus, in all
cases, except for these four, $S$ has a unique $\F$-form (in the
$SKO(n,n+1;\beta)$ case we have to assume that $\beta\in\F$ in order for
such a form to exist).

Since $H^1(Gal, O_n(\bar{\F}))$ is in canonical bijective correspondence
with classes of non-degenerate bilinear forms in $n$ variables over $\F$
(see, e.g., \cite{Serre2}), we find that
all $\F$-forms of $H(m,n)$ and $K(m,n)$ are defined by the action
on supersymplectic and supercontact forms over $\F$, respectively. In the
cases $S=S(1,2)$ and $E(1,6)$, the answer is more interesting. We construct
all $\F$-forms of these Lie superalgebras, using the theory of Lie conformal
superalgebras.

The present paper is a continuation of \cite{CantaK}, which we
refer to for terminology not explained here. The base field,
unlike in \cite{CantaK}, is an arbitrary field $\F$ of characteristic
0, and we denote by $\bar{\F}$ its algebraic closure.

In the Lie algebra case the problems considered in the present
paper were solved by Rudakov \cite{R}, whose methods we use. 

\section{$\boldsymbol{\Z}$-Gradings}\label{gradings}
In papers \cite{CantaK} and \cite{K} the base field is $\C$. However,
it is not difficult to extend all the results there to the case of an 
arbitrary algebraically closed field $\bar{\F}$ of characteristic
zero. In order to do this one has to replace exponentiable
derivations of a linearly compact algebra $S$ in the sense
of \cite{CantaK}, \cite{K}, by exponentiable derivations
in the sense of \cite{G2} (a derivation $d$ of a Lie superalgebra
$S$ over a field $\F$ is called {\em exponentiable} in the sense of
\cite{G2} if $d(H)\subset H$ for any closed $Aut S$-invariant
subspace $H$ of $S$). Also, we define the {\em group of inner automorphisms}
of $S$ to be the group generated by all elements $\exp(ad ~a)$, 
where $\exp(ad ~a)$ converges in linearly compact topology.
Then Theorem 1.7 of \cite{CantaK} on conjugacy of maximal tori
in an artinian semisimple linearly compact superalgebra still holds
over $\bar{\F}$. Consequently, the classification given in \cite{CantaK}
of primitive pairs $(L,L_0)$ up to conjugacy by inner automorphisms of
$Der L$ stands as well over $\bar{\F}$.

We first recall from \cite{CantaK} and \cite{K} the necessary
information on $\Z$-gradings of Lie superalgebras in question over
the field $\bar{\F}$. For information on finite-dimensional Lie superalgebras
we refer to \cite{K2} or \cite{K}.

Recall that $W(m,n)$ is the Lie superalgebra of all continuous derivations of 
the commutative associative superalgebra $\Lambda(m,n)=\Lambda(n)
[[x_1,\dots,x_m]]$, where $\Lambda(n)$ is the Grassmann superalgebra in $n$ odd
indeterminates $\xi_1,\dots,\xi_n$, and $x_1,\dots,x_m$ are even indeterminates.
Recall that
a $\Z$-grading of the Lie superalgebra $W(m,n)$ is called the grading of type 
$(a_1,\dots, a_m|b_1,\dots, b_n)$ if $a_i=\deg x_i=-\deg\frac{\partial}{\partial x_i}\in\N$ and
$b_i=\deg\xi_i=-\deg\frac{\partial}{\partial \xi_i}\in\Z$ 
(cf.\ \cite[Example 4.1]{K}). Every such a  grading always induces a grading 
 on the Lie superalgebra $S(m,n)$ and it induces a
grading on  $S=H(m,n)$, $K(m,n)$, $HO(n,n)$, $SHO(n,n)$,
$KO(n,n+1)$, or $SKO(n,n+1;\beta)$ if 
the defining differential
 form of $S$ is homogeneous with respect to this grading.
The induced grading on $S$ is also called a grading of type 
$(a_1,\dots, a_m|b_1,\dots, b_n)$.

The $\Z$-grading of type $(1,\dots,1|1,\dots, 1)$ is an irreducible
grading of $W(m,n)$ called its
{\em principal} grading.
In this grading $W(m,n)=\prod_{j\geq -1}\g_j$ has  $0$-th graded component isomorphic to the Lie
superalgebra $gl(m,n)$ and $-1$-st graded component isomorphic to the
standard $gl(m,n)$-module $\bar{\F}^{m|n}$. 
The even part of $\g_0$ is isomorphic to the Lie algebra $gl_m\oplus gl_n$ where
$gl_m$ (resp.\ $gl_n$) acts  trivially on $\bar{\F}^n$ (resp.\ $\bar{\F}^m$)  
and acts as the standard representation on $\bar{\F}^m$ (resp.\ $\bar{\F}^n$). 

The principal grading of $W(m,n)$ induces on $S(m,n)$,
$H(m,n)$, $HO(n,n)$ and $SHO(n,n)$,
irreducible gradings also called {\em principal}.

The $0$-th graded component of $S(m,n)$ in its principal grading
is isomorphic to the Lie superalgebra $sl(m,n)$ and its
$-1$-st graded component is isomorphic to the standard $sl(m,n)$-module
$\bar{\F}^{m|n}$. The even part of $\g_0$ is isomorphic to the Lie algebra 
$sl_m\oplus sl_n \oplus \bar{\F} c$ where
$sl_m$ (resp.\ $sl_n$) acts  trivially on $\bar{\F}^n$ (resp.\ $\bar{\F}^m$)  
and acts as the standard representation on $\bar{\F}^m$ (resp.\ $\bar{\F}^n$). 
Here $c$ acts by multiplication by $-n$ (resp.\ $-m$) on $\bar{\F}^m$
(resp.\ $\bar{\F}^n$).

Let $S=H(2k,n)=\prod_{j\geq -1}\g_j$ with its principal grading.
Recall that the Lie superalgebra $H(2k,n)$ can be identified
with $\Lambda(2k,n)/\bar{\F} 1$, where we have $2k$ even indeterminates
$q_1, \dots, q_k$, $p_1, \dots, p_k$, and $n$ odd indeterminates
$\xi_1, \dots, \xi_n$, with bracket
$[f,g]=\sum_{i=1}^k(\frac{\partial f}{\partial p_i}\frac{\partial
g}{\partial q_i}-\frac{\partial f}{\partial q_i}\frac{\partial
g}{\partial p_i})-(-1)^{p(f)}\sum_{i=1}^n\frac{\partial f}{\partial \xi_i}
\frac{\partial g}{\partial \xi_{n-i+1}}.$
 Then
$\g_0\cong spo(2k,n)$, and $\g_{-1}$ is 
isomorphic to the standard $spo(2k,n)$-module $\bar{\F}^{2k|n}$.
Here $(\g_0)_{\bar{0}}$ is spanned by elements
$\{p_ip_j, p_iq_j, q_iq_j\}$ for $i,j=1, \dots,k$, and $\{\xi_i\xi_j\}_{i\neq j}$
for $i, j=1, \dots, n$, hence it is isomorphic to $sp_{2k}\oplus so_n$.
The odd part of $\g_0$ is spanned by vectors $\{p_i\xi_j, q_i\xi_j\}$
for $i=1, \dots, k$ and $j=1, \dots, n$, hence it is isomorphic to the
$sp_{2k}\oplus so_n$-module $\bar{\F}^{2k}\otimes\bar{\F}^n$, where $\bar{\F}^{2k}$ and
$\bar{\F}^n$ are the standard $sp_{2k}$ and $so_n$-modules, respectively.
Besides, $\g_{-1}=\langle p_i, q_i, \xi_j ~|~ i=1, \dots, k, j=1, \dots, n
\rangle$, hence
$sp_{2k}$ acts trivially on $\bar{\F}^n=\langle \xi_j ~|~ j=1, \dots, n\rangle$ 
and acts as the standard
representation on $\bar{\F}^{2k}=\langle p_i, q_i ~|~ i=1, \dots,k\rangle$, and $so_n$ acts trivially
on $\bar{\F}^{2k}$ and by the standard action on $\bar{\F}^n$.

The grading of 
type $(2,1,\dots,1|1,\dots,1)$ of $W(2k+1,n)$ induces 
an irreducible grading $K(2k+1,n)=\prod_{j\geq -2}\g_j$, called the
{\em principal} grading of $K(2k+1,n)$.
Its $0$-th graded component $\g_0$ is isomorphic
to the Lie superalgebra $cspo(2k,n)$ 
and $\g_{-1}$ is isomorphic to the standard $cspo(2k,n)$-module $\bar{\F}^{2k|n}$. 

Consider the Lie superalgebra  $S=HO(n,n)=\prod_{j\geq -1}\g_j$ 
with its principal
grading.
Then $\g_0$
is isomorphic to the Lie superalgebra
$\tilde{P}(n)=\tilde{P}(n)_{-1}+\tilde{P}(n)_{0}+\tilde{P}(n)_1$, 
where $\tilde{P}(n)_{0}\cong gl_n$, and, as $gl_n$-modules,
$\tilde{P}(n)_{-1}\cong\Lambda^2(\bar{\F}^{n*})$, 
$\tilde{P}(n)_1\cong S^2(\bar{\F}^n)$, and $\g_{-1} \cong \bar
{\F}^n\oplus \bar{\F}^{n*}$, where $\bar{\F}^n$ is the standard
$gl_n$-module, and the $gl_n$-submodules $\bar{\F}^n$ and 
$\bar{\F}^{n*}$ of $\g_{-1}$ have different parities. 

The $0$-th graded component of $SHO(n,n)$ in its
principal grading is
isomorphic to the graded subalgebra $P(n)=P(n)_{-1}+P(n)_{0}+P(n)_1$
 of $\tilde{P}(n)$,  where
$P(n)_{0}\cong sl_n$, 
${P}(n)_{-1}\cong\Lambda^2(\bar{\F}^{n*})$ and 
${P}(n)_1\cong S^2(\bar{\F}^n)$,
 and its $-1$-st graded component is
isomorphic to the standard $P(n)$-module $\bar{\F}^n\oplus \bar{\F}^{n*}$.

The $\Z$-grading of type $(1, \dots, 1|1, \dots,1,2)$ of $W(n,n+1)$
induces on the Lie superalgebras $KO(n,n+1)$ and
$SKO(n,n+1;\beta)$ an irreducible grading
called {\em principal}.
In these cases the $\g_0$-module $\g_{-1}$ is obtained from that
of $SHO(n,n)$ by adding some operators which act as scalars on
$\bar{\F}^n$ and $\bar{\F}^{n*}$.

The $\Z$-grading of type $(1, \dots,1|0,\dots,0)$
of $S=W(m,n)$, $S(m,n)$, $HO(n,n)$,
$SHO(n,n)$, the $\Z$-gradings of type  $(1,\dots,1|2, \dots,2,0,\dots,0)$
and $(2,1,\dots,1|$ $2, \dots,2,$ $0,\dots,0)$ with $h$ zeros
of $S=H(m,2h)$ and $K(m,2h)$, respectively, and the
$\Z$-grading
of type $(1, \dots,1|0,\dots,0,1)$ of $S=KO(n,n+1)$,
$SKO(n,n+1;\beta)$,
is called the {\em subprincipal} grading of $S$.

The Lie superalgebra $S=SKO(2,3;1)=\prod_{j\geq -1}\g_j$ in its subprincipal grading
has $0$-th graded component $\g_0$ isomorphic to the
semidirect sum of $S(0,2)$ and the subspace of $\Lambda(2)$ spanned
by all the monomials except for $\xi_1\xi_2$,
and $\g_{-1}\cong \Lambda(2)$. The even
part of $\g_0$ is isomorphic to $sl_2 \oplus \bar{\F}$ and, 
as an $sl_2$-module,
$\g_{-1}=\bar{\F}^2\oplus \bar{\F}^2$, where the two copies of $\bar{\F}^2$ have different
parities and $sl_2$ acts by the standard action on the even copy and
 trivially on the odd copy. The algebra  of outer
derivations of $S$ is isomorphic to $sl_2$ (cf.\ \cite[Remark 4.15]{CantaK}); 
it acts trivially on
the even subspace of $\g_{-1}$ and by the standard action on the
odd one. Finally, $\bar{\F}$ acts on $\g_{-1}$ by multiplication by $-2$.

The $\Z$-grading of $W(1,2)$ of type $(2|1,1)$ induces
a grading on $S(1,2)=\prod_{j\geq -2}\g_j$, which is not irreducible.
Then $\g_0\cong sl_2 \oplus \bar{\F} c$
where $c$ acts on $S$ as the grading operator, and
$\g_{-1}=\bar{\F}^2\oplus \bar{\F}^2$, where $\bar{\F}^2$ is the standard $sl_2$-module.
The two copies of the standard $sl_2$-module in $\g_{-1}$ are both
odd.

Likewise, the $\Z$-grading of $W(3,3)$ of type $(2,2,2|1,1,1)$
induces a grading on $SHO(3,3)=\prod_{j\geq -2}\g_j$ , 
which is not irreducible.
Here $\g_0\cong sl_3$ and $\g_{-1}=\bar{\F}^3\oplus\bar{\F}^3$, where
$\bar{\F}^3$ is the standard $sl_3$-module. The two copies of the 
standard $sl_3$-module in $\g_{-1}$ are both odd.

Consider the Lie superalgebra
$K(1,6)=\prod_{j\geq -2}\mathfrak{g}_j$
with its principal grading. Then
$\g_0=sl_4\oplus\bar{\F} c$ and 
$\mathfrak{g}_{-1}\cong\Lambda^2\bar{\F}^4$, where 
$\bar{\F}^4$ denotes the standard $sl_4$-module, 
$\g_1\cong\mathfrak{g}_{-1}^*\oplus\mathfrak{g}_1^+\oplus \mathfrak{g}_1^-$, as
$sl_4$-modules, with
$\mathfrak{g}_1^+\cong S^2\bar{\F}^4$ and $\mathfrak{g}_1^-\cong S^2(\bar{\F}^{4^*})$.
The Lie superalgebra $S=E(1,6)$ is the graded subalgebra of $K(1,6)$ generated by
$\mathfrak{g}_{-1}+\mathfrak{g}_0+(\mathfrak{g}_{-1}^*+\mathfrak{g}_1^+)$
 (cf.\ \cite[Example 5.2]{K}, \cite[\S 4.2]{CK}, \cite[\S 3]{S}).
It follows that the $\Z$-grading of type $(2|1,1,1,1,1,1)$ induces on
$E(1,6)$ an irreducible grading, called the {\em principal} grading
of $E(1,6)$,
where
 ${\mathfrak{g}}_{-1}=\langle \xi_i, \eta_i
\rangle$,
 $\mathfrak{g}_{-1}^*=\langle t\xi_i, t\eta_i\rangle$ and 
$\mathfrak{g}_1^+=\langle \xi_1\xi_2\xi_3, \xi_1\eta_2\eta_3,
\xi_2\eta_1\eta_3, \xi_3\eta_1\eta_2, \xi_1(\xi_2\eta_2+\xi_3\eta_3),
\xi_2(\xi_1\eta_1+\xi_3\eta_3), \eta_3(\xi_1\eta_1-\xi_2\eta_2),
\xi_3(\xi_1\eta_1+\xi_2\eta_2), \eta_2(\xi_1\eta_1-\xi_3\eta_3),
\eta_1(\xi_2\eta_2-\xi_3\eta_3)\rangle$, and
$\mathfrak{g}_1^-$ is obtained from $\mathfrak{g}_1^+$ by exchanging
$\xi_i$ with $\eta_i$ for every $i=1,2,3$.

Next, the
{\em principal}
grading of $E(3,6)$ is an irreducible grading of depth two whose
$0$-th graded component is isomorphic to
$sl_3\oplus sl_2\oplus\mathbb{C}c$, and 
whose $-1$-st graded component is
isomorphic, as an  $sl_3\oplus sl_2$-module, to $\bar{\F}^3\boxtimes \bar{\F}^2$ 
where $\bar{\F}^3$ and $\bar{\F}^2$ denote the standard $sl_3$ and $sl_2$-modules,
respectively. Here $c$ acts on $E(3,6)$ as the grading operator
(with respect to its principal grading).
Likewise, the {\em principal}
grading of $E(3,8)$ is an irreducible grading of depth three whose
$0$-th graded component is isomorphic to
$sl_3\oplus sl_2\oplus\mathbb{C}c$ , and 
whose $-1$-st graded component is
isomorphic, as an  $sl_3\oplus sl_2$-module, to $\bar{\F}^3\boxtimes \bar{\F}^2$ 
where $\bar{\F}^3$ and $\bar{\F}^2$ denote the standard $sl_3$ and $sl_2$-modules,
respectively; $c$ acts on $E(3,8)$ as the grading operator.

The Lie superalgebra $S=E(4,4)$ has even part isomorphic to $W_4$
and odd part isomorphic to the $W_4$-module $\Omega^1(4)^{-\frac{1}{2}}$.
The bracket between two odd elements $\omega_1$ and $\omega_2$
is defined as: $[\omega_1, \omega_2]=d\omega_1\wedge \omega_2+
\omega_1\wedge d\omega_2$.
The {\em principal} grading of $S$ is an 
irreducible $\Z$-grading of depth 1 whose $0$-th graded
component $\g_0=\langle x_i\frac{\partial}{\partial x_j}, x_idx_j\rangle$
is isomorphic to the Lie superalgebra $\hat{P}(4)$
and 
$\g_{-1}=\langle \frac{\partial}{\partial x_j}, dx_j\rangle$
is isomorphic to 
the standard $\hat{P}(4)$-module $\bar{\F}^{4|4}$.
We recall that $\hat{P}(4)=P(4)+\bar{\F} z$ is a (non-trivial)
central extension of $P(4)$ with center $\bar{\F} z$ (see \cite{CantaK}, \cite{K},
\cite{S}).


Finally, the {\em principal} grading of the Lie superalgebra $E(5,10)$ 
is irreducible of depth 2,
with $0$-th graded component isomorphic to
$sl_5$ and $-1$-st graded component isomorphic ro
 $\Lambda^2 \bar{\F}^5$, where $\bar{\F}^5$
is the standard $sl_5$-module.

\medskip

Given a simple infinite-dimensional linearly compact  Lie superalgebra 
$S=\prod_{j\geq -d}\g_j$
with its principal or subprincipal grading, we will call $S_0=\prod_{j\geq 0}
\g_j$ the principal or subprincipal subalgebra of $S$, 
respectively. Likewise,
if $S=\prod_{j\geq -d}\g_j$
with a grading of a given type, we will call $S_0=\prod_{j\geq 0}
\g_j$ the subalgebra of $S$ of this  type.

\begin{remark}\label{newconstruction}\em One can show that 
every non-graded maximal open subalgebra of
any non-exceptional simple infinite-dimensional linearly compact
Lie superalgebra $S$ in its defining embedding in $W(m,n)$, can be
constructed as the intersection of $S$ with a graded subalgebra
of $W(m,n)$. For example, the maximal open subalgebra
$L_0(0)$ of $S=H(m,1)$ constructed in \cite[Example 3.3]{CantaK},
is the intersection of $S$ with the subprincipal subalgebra
of $W(m,1)$. Since the supersymplectic form is not homogeneous
with respect to the subprincipal grading of $W(m,1)$, $L_0(0)$ is not
graded. We shall call this subalgebra
 the {\em subprincipal} subalgebra of $H(m,1)$.
\end{remark}

If $\mathfrak{g}$ is a Lie algebra acting linearly on a vector
space $V$ over $\bar{\F}$, we denote by $\exp(\g)$ the linear
algebraic subgroup of $GL(V)$, generated by all $\exp a$, where 
$a$ is a (locally) nilpotent endomorphism of $V$, and by $t^a$, 
where $a$ is a diagonalizable endomorphism of $V$ with integer eigenvalues and $t\in\bar{\F}^\times$. 

If a group $G$ is an almost direct product of two subgroups $G_1$ and
$G_2$ (i.e., both $G_1$ and $G_2$ are normal subgroups
and $G_1\cap G_2$ is a finite central subgroup of $G$)
 we will denote it by $G=G_1\cdot G_2$.
We will often make use of the following simple result:
\begin{proposition}\label{maximal} Suppose we have a representation of a Lie
superalgebra $\g$ over $\bar{\F}$ in a vector superspace $V$, and a
faithful representation of
a group $G$ in $V$, containing $\exp(\g_{\bar{0}})$, preserving parity and such
that conjugation by elements of $G$ induces automorphisms of $\g$. Then
the maximal possible $G$ are as follows in the following cases: 
\begin{itemize}
\item[$(a)$] if $\g=sl_n$ and $V=\bar{\F}^n\oplus \bar{\F}^n$
with the same parity, then $G$ is an almost direct product of
$GL_2$ and $SL_n$; in particular if $n=2$ then $G=\bar{\F}^\times\cdot SO_4$
and if $n=3$ then $G=\bar{\F}^\times\cdot(SL_2\times SL_3)$; 
\item[$(b)$] if $\g=sl_2\oplus sl_3$ and $V=\bar{\F}^3\boxtimes \bar{\F}^2$,
then $G=\bar{\F}^\times\cdot(SL_2\times SL_3)$;
\item[$(c)$] if $\g=sl_5$ and $V=\Lambda^2\bar{\F}^5$, then $G=GL_5$;
\item[$(d)$] if $\g=sl(m,n)$ and $V=\bar{\F}^{m|n}$ is the standard $sl(m,n)$-module, then $G=GL_m\times GL_n$;
\item[$(e)$] if $\g=spo(2k,n)$ and $V=\bar{\F}^{2k|n}$ is the standard $spo(2k,n)$-module, then 
$G=\bar{\F}^\times\cdot(Sp_{2k}\times O_n)$;
\item[$(f)$] if $\g=P_n$ and $V=\bar{\F}^n+\bar{\F}^{n*}$ is the standard $P_n$-module, 
then $G=\bar{\F}^\times\cdot GL_n$; 
\item[$(g)$] if $\g=\hat{P}_4$ and $V=\bar{\F}^4+\bar{\F}^{4*}$ is the standard $\hat{P}_4$-module, then $G=GL_4$.
\end{itemize}
In all cases when $G=\bar{\F}^\times\cdot G_1$, the group $\bar{\F}^\times$ acts on $V$
by scalar multiplication.
\end{proposition}
{\bf Proof.} 
Consider the map $f: G\longrightarrow Aut(\g)$ that
associates to every element of $G$ the induced automorphism by conjugation of
$\g$. Then the kernel of $f$ consists of the elements of $G$ commuting with
$\g$.
Suppose, as in $(a)$, that $\g=sl_n$ and $V=\bar{\F}^n\oplus\bar{\F}^n$, where the
two copies of the standard $sl_n$-module have the same parity.
Then $Im f$ consists of inner automorphisms of $sl_n$, i.e.,
$Im f=PGL_n$, and $\ker f=GL_2$. We have therefore the following
exact sequence:
$$1\longrightarrow GL_2\longrightarrow G\longrightarrow PGL_n\longrightarrow 1.$$
Since there is in $G$ a complementary to $GL_2$ subgroup,
which is $SL_n$, we conclude that $G$ is an almost direct product of $GL_2$
and $SL_n$. It follows that if $n$ is odd, then $G=\bar{\F}^\times\cdot(SL_2\times SL_n)$, 
and if $n$ is even then $G=\bar{\F}^\times\cdot((SL_2\times SL_n)/C_2)$
where $C_2$ is the cyclic subgroup of order two of
$SL_2\times SL_n$ generated by $(-I_2, -I_n)$, proving $(a)$.
The same argument proves $(b)$.

By the same argument, in case $(c)$ we get the exact sequence
$$1\longrightarrow \bar{\F}^\times \longrightarrow G\longrightarrow PGL_5\longrightarrow 1.$$
Since $G$ contains a complementary to $\bar{\F}^\times$ subgroup, which
is $SL_5$,
 we conclude that $G=GL_5$. 

If $\g=sl(m,n)$ and $V=\bar{\F}^{m|n}$ is its standard representation,
then $gl_m$ acts irreducibly on $\bar{\F}^m$ and
$gl_n$ acts irreducibly on $\bar{\F}^n$, hence $G=GL_m\times GL_n$, proving $(d)$.

Suppose that $\g=spo(2k,n)$ and $V=\bar{\F}^{2k|n}$ is the standard 
$spo(2k,n)$-module. Define on $\g_{-1}=\langle p_i, q_i, \xi_j~|~i=1, \dots, k,
 ~j=1, \dots, n\rangle$ the following symmetric
bilinear form: $(p_i, q_j)=\delta_{i,j}$, $(\xi_i, \xi_j)=\delta_{i, n-j+1}$,
$(p_i, p_j)=0=(q_i, q_j)=(p_i, \xi_j)=(q_i, \xi_j)$. Then
$G$ consists of the automorphisms of $\bar{\F}^{2k}\oplus\bar{\F}^n$ preserving
the bilinear form $(\cdot, \cdot)$ up to multiplication by a scalar,
hence $G=\bar{\F}^\times\cdot(Sp_{2k}\times O_n)$, proving $(e)$.

Finally, let $\g=P_n$ and $V=\bar{\F}^n+\bar{\F}^{n*}$ be the standard $P_n$-module.
Then $gl_n$ acts irreducibly on $\bar{\F}^n$ and $\bar{\F}^{n*}$ which have
different parities. It follows that $G=\bar{\F}^\times\cdot GL_n$. Likewise,
if $\g=\hat{P}(4)$, then $G=\bar{\F}^\times\cdot SL_4$ since the group of automorphisms
of $\g$ is $SL_4$.
\hfill$\Box$
\section{On $\boldsymbol{Aut S}$}\label{general}
Let $S$ be a linearly compact infinite-dimensional Lie superalgebra
over $\bar{\F}$ and let $Aut S$ denote the group of all continous automorphisms
of $S$.
Let  $S=S_{-d}\supset\dots
\supset S_0\supset \dots$ be a filtration of $S$ by open subalgebras
 such that all
$S_j$ are $Aut S$-invariant and $Gr S=\oplus_{j\geq -d}\g_j$ is a transitive 
graded Lie superalgebra.
Denote by $Aut(Gr S)$ the group of automorphisms of $Gr S$ preserving
the grading. Denote by $Autf S$ the subgroup consisting of $g\in Aut S$
 which induce
an identity automorphism of $Gr S$, and let $Autgr S$ be the subgroup of 
$Aut(Gr S)$ consisting of automorphisms induced by $g\in Aut S$. We have
an exact sequence:
\begin{equation}
1\rightarrow Autf S\rightarrow Aut S\rightarrow Autgr S\rightarrow 1.
\label{exact}
\end{equation}
\begin{proposition}\label{Victor}
(a) The restriction map $Autgr S\rightarrow GL(\g_{-1})$ is injective.

\noindent
(b) $Autf S$ consists of inner automorphisms of $Der S$. In fact 
$Autf S=\exp ad$ $(Der S)_1$, where $(Der S)_1$ is the first member
of the filtration of $Der S$, induced by that of $S$.

\noindent
(c) If $S_0$ is a graded subalgebra, i.e., $S_j=\g_j \oplus S_{j+1}$ for all $j \geq -d$ such that
$[\g_i,\g_j]\subset\g_{i+j}$, then $Autgr S = Aut(Gr S)$ and
$$Aut S=Autf S\rtimes Aut(Gr S).$$
\end{proposition}
{\bf Proof.} By transitivity, $\g_{-n}=\g_{-1}^n$ for $n \geq 1$, and, in 
addition, we have 
the well-known injective $Autgr S$-equivariant
map $\g_n\rightarrow Hom(\g_{-1}^{\otimes(n+1)}, \g_{-1})$ for $n \geq 0$, which
implies $(a)$.

If $\sigma\in Autf S$, then $\sigma=1+\sigma_1$, where $\sigma_1(L_j)\subset
L_{j+1}$. Hence $\log\sigma=\sum_{n\geq 1}(-1)^{n+1}\frac{\sigma_1^n}{n}$
converges and $e^{t\log\sigma}$ converges to a one-parameter
subgroup of $Autf S$. Hence $\sigma$ is an inner automorphism of $Der S$,
proving $(b)$.

If $S_0$ is a graded subalgebra, we have an obvious inclusion $Aut(Gr S)
\subset Aut S$ and exact sequence (\ref{exact}), proving $(c)$.
\hfill$\Box$
\begin{remark}\em  By Proposition \ref{Victor}(a), $Autgr S$ is a subgroup of $GL(\g_{-1})$
 whose Lie algebra is $Gr_0 Der S$ acting on $\g_{-1}$. By Proposition \ref{Victor}(b),
$Autf S$ is a prounipotent group.
\end{remark}


\section{Invariant Subalgebras}\label{invsub}
Given a linearly compact Lie superalgebra $L$, we call
{\em invariant}
a subalgebra of $L$ which is 
invariant with respect to all its inner automorphisms,
or, equivalently, which  contains
all elements
$a$ of $L$ such that $\exp(ad(a))$ converges in the linearly compact
topology. 
It turns out that an open
subalgebra of minimal codimension in a
linearly compact infinite-dimensional simple Lie superalgebra $S$
over $\bar{\F}$ is always invariant
under all inner automorphisms of $S$ (see \cite{CantaK}).
\begin{example}\label{SHOtilde}\em  
We recall that the Lie superalgebra $S=SHO^\sim(n,n)$
is the subalgebra of $HO(n,n)$ defined as follows:
$$SHO^\sim(n,n)=\{X\in HO(n,n)~|~ X(F v)=0\}$$
where $v$ is the volume form associated to the usual
 divergence and $F=1-2\xi_1\dots\xi_n$ (cf.\ \cite[\S 5]{CantaK}).
Let $S_0$ be the intersection of $S$ with the principal
subalgebra of $W(n,n)$. Then the Weisfeiler filtration
associated to $S_0$ has depth one and
$\overline{Gr S}\cong SHO^{\prime}(n,n)$  with  the
  $\Z$-grading of type  $(1,\dots, 1|1, \dots,1)$ 
(cf.\  \cite[Example 5.2]{CantaK}).
Here and further by $\overline{Gr S}$ we denote the completion of the graded
Lie superalgebra associated to the above filtration.
By \cite[Proposition 1.11]{CantaK}, $S_0$ is a maximal open
subalgebra of $S$. It is easy to see that it is also
an invariant subalgebra. This subalgebra is called
the {\em principal} subalgebra of $S$.
\end{example}
\begin{example}\label{SKOtilde}\em 
We recall that the Lie superalgebra $S=SKO^\sim(n,n+1)$ is the
subalgebra of $KO(n,n+1)$ defined as follows:
$$SKO^\sim(n,n+1)=\{X\in KO(n,n+1)~|~X(F v_{\beta})=0\}$$
where $v_{\beta}$ is the volume form attached to
the divergence $div_{\beta}$ for $\beta=(n+2)/n$
 and $F=1+\xi_1\dots\xi_{n+1}$ (cf.\ \cite[\S 5]{CantaK}).
Let $S_0$ be the intersection of $S$ with the subalgebra of $W(n,n+1)$
of type $(1,\dots,1|1,\dots,1,2)$. Then the Weisfeiler
filtration associated to $S_0$ has depth 2 and
$\overline{Gr
  S}\cong SKO(n,n+1; (n+2)/n)$ with  its
  principal grading.
By \cite[Proposition 1.11]{CantaK},
$S_0$ is a maximal open subalgebra of $S$. It is easy to see that 
it is also an invariant subalgebra. This subalgebra is called
the {\em principal} subalgebra of $S$.
\end{example}

 A complete list of
invariant maximal open subalgebras  in
all simple linearly compact infinite-dimensional
Lie superalgebras over $\bar{\F}$, is given in the following theorem
(cf.\ \cite[Theorem 11.1]{CantaK}):

\begin{theorem}\label{invariant}
The following is a complete list of invariant maximal open subalgebras
in infinite-dimensional linearly compact simple Lie superalgebras $S$
over $\bar{\F}$:
\begin{description}
\item[$(a)$] the principal subalgebra of $S$;
\item[$(b)$] the subprincipal subalgebra of $S=W(m,1)$,
$S(m,1)$, $H(m,1)$, $H(m,2)$, $K(m,2)$, $KO(2,3)$, $SKO(2,3;\beta)$,
the subalgebra of type $(2,1,\dots,1|0,2)$ of $K(m,2)$
and the subalgebra of type $(1,\dots,1|0,2)$ of $H(m,2)$;
\item[$(c)$] the subalgebra of type $(1,1|-1,-1,0)$ of
$SKO(2,3;\beta)$ for $\beta\neq 1$;
\item[$(d)$] the subalgebras of $S=S(1,2)$, $SHO(3,3)$,
$SKO(2,3;1)$ conjugate to the principal subalgebra 
 by the subgroup of $Aut S$ generated by the automorphisms
$\exp(ad(\mathfrak{a}))$ where $\mathfrak{a}$ is the algebra of
 outer derivations of $S$;
\item[$(e)$] the subalgebras of $S=SKO(3,4;1/3)$ conjugate to the
subprincipal subalgebra  by the  automorphisms
$\exp(ad(t\xi_1\xi_2\xi_3))$ with $t\in\bar{\F}$.
\end{description}
\end{theorem}

\begin{theorem}\label{autSinv} The following is a complete list of
maximal among the open $Aut S$-invariant
subalgebras in
infinite-dimensional linearly compact simple Lie superalgebras $S$:
\begin{description}
\item[$(a)$] the principal subalgebra of
$S\neq S(1,2)$, $SHO(3,3)$, 
and $SKO(2,3;1)$;
\item[$(b)$]  the subprincipal subalgebra 
of $S=W(m,1)$, $S(m,1)$, $H(m,1)$, $KO(2,3)$, and
$SKO(2,3;\beta)$;
\item[$(c)$] the subalgebra of type $(1,1|-1,-1,0)$ in $S=SKO(2,3;\beta)$
with $\beta\neq 1$;
\item[$(d)$] the subalgebras of type 
$(2|1,1)$ and $(2,2,2|1,1,1)$ in 
$S=S(1,2)$ and\break $SHO(3,3)$, 
respectively.
\end{description}
\end{theorem}
{\bf Proof.} 
We will prove  that the subalgebras  listed in $(a)-(d)$
are $Aut S$-invariant, and,
in order to show that they
exhaust all maximal among open $Aut S$-invariant subalgebras
of $S$, it suffices to show that for every subalgebra $S'_0$ of $S$
listed in Theorem \ref{invariant},
 $\cap_{\varphi\in Aut S}\varphi(S'_0)$
is contained in one of them.
Indeed, if $S_0$ is a maximal among the $Aut S$-invariant
open subalgebras of $S$, then $S_0$ is an invariant
subalgebra of $S$, hence, every maximal open subalgebra $S'_0$
of $S$ containing $S_0$, is  invariant, i.e., $S'_0$ is
one of the subalgebras of $S$ listed in Theorem \ref{invariant}.
Therefore, $S_0\subset \cap_{\varphi\in Aut S}\varphi(S'_0)$.

If $S\neq W(m,1)$, $S(m,1)$, $H(m,1)$, $H(m,2)$,
$K(m,2)$, $KO(2,3)$, $SKO(2,3;\beta)$, 
  $S(1,2)$, $SHO(3,3)$ and $SKO(3,4;1/3)$,
 then, in view of Theorem
\ref{invariant}, 
the principal subalgebra of $S$ 
is the unique invariant maximal open subalgebra 
of $S$, hence it is
invariant with respect to all automorphisms of $S$
and it is the unique maximal among open $Aut S$-invariant
subalgebras of $S$.

If $S=W(m,1)$, $S(m,1)$, or $KO(2,3)$, then, 
according to Theorem \ref{invariant}, $S$ has two invariant subalgebras: the
 principal and subprincipal subalgebras. 
These two subalgebras have different
codimension hence, each of them is invariant 
with respect to all automorphisms of $S$. 

If $S=H(m,1)$, then $S$ has two invariant maximal open subalgebras: the
principal and the subprincipal subalgebras.
Since the principal subalgebra
is graded and the subprincipal subalgebra is not 
(see Remark \ref{newconstruction}),
each of them is invariant with respect to all
automorphisms of $S$.

If $S=H(m,2)$, then the principal subalgebra 
is the unique subalgebra of $S$ of minimal codimension, hence it
is $Aut S$-invariant. Besides, it
is the unique maximal among $Aut S$-invariant
subalgebras of $S$. Indeed, the invariant subalgebras of $S$
of type $(1, \dots,1|2,0)$ and  $(1, \dots,1|0,2)$ are conjugate
by an outer automorphism of $S$ and their intersection
is contained in the principal subalgebra of $S$.
By the same arguments, if $S=K(m,2)$, then the principal subalgebra
of $S$ is its unique maximal among open $Aut S$-invariant
subalgebras.

If $S=SKO(3,4;1/3)$, then, according to
Theorem \ref{invariant}, $S$ has infinitely many invariant subalgebras
which are conjugate  to the subprincipal 
subalgebra. Besides, the principal subalgebra of $S$ 
 is an invariant subalgebra. Note that the principal grading of $S$
has depth 2 and the subprincipal grading of $S$ 
has depth 1, therefore the 
principal and subprincipal subalgebras are not conjugate.
It follows that the principal subalgebra
is invariant with respect to all
automorphisms of $S$. 
In fact, it
 is the unique maximal among $Aut S$-invariant subalgebras of $S$,
since the intersection of all the subalgebras of $S$ listed in
Theorem \ref{invariant}$(e)$ is the subalgebra of $S$ of type $(2,2,2|1,1,1,3)$,
which is contained in the principal subalgebra.

If $S=SKO(2,3;1)$, then $S$ has infinitely many invariant subalgebras
which are conjugate to the principal subalgebra, besides, the
subprincipal subalgebra is also an invariant subalgebra of $S$.
The  principal and subprincipal subalgebras
have codimension $(2|3)$ and $(2|2)$, respectively,
hence they cannot be conjugate. It follows that
the subprincipal subalgebra is invariant with respect to
all automorphisms of $S$. 
In fact, it is the unique maximal among
$Aut S$-invariant subalgebras of $S$, since it contains the 
intersection of all subalgebras which are conjugate
to the principal subalgebra
(cf.\ \cite[Remark 4.16]{CantaK}).

If $S=SKO(2,3;\beta)$ for $\beta\neq 0,1$, then,
according to Theorem \ref{invariant}, $S$ has three invariant
maximal open subalgebras, i.e., the subalgebras
of type $(1,1|1,1,2)$, $(1,1|0,0,1)$ and $(1,1|-1,-1,0)$.
The subalgebras of type $(1,1|1,1,2)$ and $(1,1|-1,-1,0)$
have codimension $(2|3)$ and the subalgebra of type
$(1,1|0,0,1)$ has codimension $(2|2)$. It follows that
the subprincipal subalgebra is invariant with respect to
all automorphisms of $S$, since it is the unique 
subalgebra of  minimal codimension.
Consider the grading of $S$ of type $(1,1|1,1,2)$:
this is an irreducible grading of depth 2, whose $0$-th
graded component is isomorphic to the Lie
superalgebra $\tilde{P}(2)=P(2)+\bar{\F}(\xi_3+\beta\Phi)$. Its
$-2$-nd graded component is $\bar{\F} 1$, on which $\xi_3+\beta\Phi$ acts
as the scalar $-2$. Its $-1$-st graded component is spanned by
vectors $\{x_i\}$ and $\{\xi_i\}$, with $i=1,2$, hence it
is isomorphic to the standard
$P(2)$-module, and $\xi_3+\beta\Phi$ acts on $\sum_{i=1}^2\bar{\F} x_i$
(resp.\ $\sum_{i=1}^2\bar{\F} \xi_i$) as the scalar $-1+\beta$
(resp.\ $-1-\beta$).
Now let us consider the grading of $S$ of type $(1,1|-1,-1,0)$:
this is an irreducible grading of depth 2, whose $0$-th
graded component is isomorphic to the Lie
superalgebra $\tilde{P}(2)=P(2)+\bar{\F}(\xi_3+\beta\Phi)$. Its
$-2$-nd graded component is $\bar{\F} \xi_1\xi_2$, on which $\xi_3+\beta\Phi$ acts
as the scalar $-2\beta$. Its $-1$-st graded component is spanned by
vectors $\{\xi_i(\xi_3+(2\beta-1)\Phi)\}$ and $\{\xi_i\}$, with $i=1,2$, hence it
is isomorphic to the standard
$P(2)$-module, and $\xi_3+\beta\Phi$ acts on $\sum_{i=1}^2\bar{\F} \xi_i(\xi_3+
(2\beta-1)\Phi)$
(resp.\ $\sum_{i=1}^2\bar{\F} \xi_i$) as the scalar $1-\beta$
(resp.\ $-1-\beta$). Since we assumed $\beta\neq 1$, the two
gradings are not isomorphic, hence the subalgebras of
type $(1,1|1,1,2)$ and $(1,1|-1,-1,0)$ are not conjugate by
any automorphism of $S$. We conclude that they are invariant
with respect to all automorphisms of $S$.

Finally, if  $S=S(1,2)$ or $S=SHO(3,3)$, then,
each of the subalgebras listed in $(d)$
is the intersection of all
invariant maximal open subalgebras of $S$, which lie
in an $Aut S$-orbit, by Theorem \ref{invariant},
thus it is the unique maximal
among $Aut S$-invariant subalgebras of $S$.
\hfill$\Box$


\section{The Group $\boldsymbol{Autgr S}$}
In this section, for every simple infinite-dimensional linearly compact Lie 
superalgebra $S$ over $\bar{\F}$, we fix the following maximal among
$Aut S$-invariant open subalgebras of $S$, which we shall denote by $S_0$:
\begin{enumerate}
\item the principal subalgebra of $S\neq S(1,2)$, $SHO(3,3)$, 
$SKO(2,3;1)$; 
\item the subalgebra of type $(2|1,1)$ in $S=S(1,2)$;
\item the subalgebra of type $(2,2,2|1,1,1)$ in $S=SHO(3,3)$;
\item the subprincipal subalgebra of $S=SKO(2,3;1)$.
\end{enumerate}

\begin{remark}\label{canonical}\em In \cite[\S 11]{CantaK} we
introduced the notion of the {\em canonical} subalgebra of
$S$,
defined as the intersection of all subalgebras of minimal codimension
in $S$. It follows from the definition that the canonical subalgebra
of $S$ is an $Aut S$-invariant subalgebra. If $S\neq KO(2,3)$, $SKO(2,3;\beta)$
with $\beta\neq 1$, and $S\neq SKO(3,4;1/3)$, then the maximal
among $Aut S$-invariant subalgebras 
$S_0$ of $S$ we have chosen is the canonical subalgebra of $S$.
\end{remark}

\medskip

Let $S_{-1}$ be a minimal subspace of $S$, properly containing
the subalgebra $S_0$
and invariant with respect to the group $Aut S$, and let
$S=S_{-d} \supsetneq
S_{-d+1} \supset \cdots \supset S_{-1} \supset S_0 \supset
\cdots$ be the associated Weisfeiler filtration of
$S$.  
All members of the Weisfeiler  filtration associated to $S_0$
 are invariant with respect
to the group $Aut S$. Let $Gr S=\oplus_{j\geq -d}\g_j$ be the associated 
$\Z$-graded Lie superalgebra. 
In this section we will describe the group
$Autgr S$ introduced in Section \ref{general}, for every $S$.
The results are summarized in Table 1 (where by
$\Pi V$ we denote $V$ with reversed parity).

\bigskip

\begin{table}[htbp]
\begin{center}
\begin{tabular}{|c|c|c|c|}
\hline
$S$ & $\g_0$ & $\g_0$-module $\g_{-1}$ & $Autgr S$ \\
\hline
$W(m,n)$, $(m,n)\neq (1,1)$ & $gl(m,n)$ & $\bar{\F}^{m|n}$ & $GL_m\times GL_n$\\
$S(1,2)$  & $sl_2\oplus\bar{\F}$ & $\Pi\bar{\F}^2\oplus\Pi\bar{\F}^2$ & $\bar{\F}^\times\cdot SO_4$ \\
$S(m,n)$, $(m,n)\neq (1,2)$ & $sl(m,n)$ & $\bar{\F}^{m|n}$  & $GL_m\times GL_n$ \\
$H(2k,n)$ & $spo(2k,n)$ & $\bar{\F}^{2k|n}$ & $\bar{\F}^\times\cdot(Sp_{2k}\times O_n)$ \\
$K(2k+1,n)$  & $cspo(2k,n)$ & $\bar{\F}^{2k|n}$ & $\bar{\F}^\times\cdot(Sp_{2k}\times O_n)$ \\
$HO(n,n)$, $n>2$ & $\tilde{P}(n)$ & $\bar{\F}^n\oplus\Pi\bar{\F}^{n*}$ & $\bar{\F}^\times\cdot GL_n$  \\
$SHO(3,3)$ & $sl_3$ & $\Pi\bar{\F}^3\oplus\Pi\bar{\F}^3$ & $\bar{\F}^\times\cdot(SL_3\times SL_2)$\\ 
$SHO(n,n)$, $n>3$  & $P(n)$ & $\bar{\F}^n\oplus\Pi\bar{\F}^{n*}$ & $\bar{\F}^\times\cdot GL_n$  \\
$KO(n,n+1)$, $n\geq 2$ & $c\tilde{P}(n)$ &  $\bar{\F}^n\oplus\Pi\bar{\F}^{n*}$ & 
$\bar{\F}^\times\cdot GL_n$  \\
$SKO(2,3;1)$ & $\langle 1, \xi_1, \xi_2\rangle\rtimes S(0,2)$ & $\Lambda(2)$
& $\bar{\F}^\times\cdot(SL_2\times SL_2)$\\
$SKO(2,3;\beta)$, $\beta\neq 0,1$ & $\tilde{P}(2)$ & $\bar{\F}^2\oplus\Pi\bar{\F}^{2*}$  & $\bar{\F}^\times\cdot GL_2$\\
$SKO(n,n+1;\beta)$, $n>2$ & $\tilde{P}(n)$ & $\bar{\F}^n\oplus\Pi\bar{\F}^{n*}$  & $\bar{\F}^\times\cdot GL_n$\\
$SHO^\sim(n,n)$, $n>2$ & $P(n)$ & $\bar{\F}^n\oplus\Pi\bar{\F}^{n*}$ & $SL_n$\\
$SKO^\sim(n,n+1)$ & $\tilde{P}(n)$ & $\bar{\F}^n\oplus\Pi\bar{\F}^{n*}$ & $SL_n$\\
 $E(1,6)$ & $so_6\oplus\bar{\F}$ & $\Pi\bar{\F}^6$ & $\bar{\F}^\times\cdot SO_6$\\
 $E(3,6)$ & $sl_3\oplus sl_2\oplus\bar{\F}$ & $\Pi(\bar{\F}^3\boxtimes \bar{\F}^2)$  & $\bar{\F}^\times\cdot(SL_2 \times SL_3)$ \\
 $E(5,10)$ & $sl_5$ & $\Pi(\Lambda^2 \bar{\F}^5)$ & $GL_5$ \\
$E(4,4)$ & $\hat{P}(4)$ & $\bar{\F}^{4|4}$ & $GL_4$ \\
 $E(3,8)$ & $sl_3\oplus sl_2\oplus\bar{\F}$ & $\Pi(\bar{\F}^3\boxtimes \bar{\F}^2)$  & $\bar{\F}^\times\cdot(SL_2 \times SL_3)$ \\
\hline
\end{tabular}
\begin{center}
\textbf{Table 1.}
\end{center}
\end{center}
\end{table}

\begin{theorem}\label{AutgrS}
Let $S$ be a simple infinite-dimensional linearly compact Lie superalgebra
over $\bar{\F}$.
Then $Autgr S$ is the 
algebraic group listed in the last column of Table 1.
In all cases when $Autgr S\cong \bar{\F}^\times\cdot G_1$, the
group $\bar{\F}^\times$ acts on $\g_{-1}$
by scalar multiplication.
\end{theorem}
{\bf Proof.}
Let $S=W(m,n)$ with $(m,n)\neq (1,1)$ or $S=S(m,n)$ with  $(m,n)\neq (1,2)$,
with the principal grading.
By Proposition  \ref{maximal}$(d)$, 
$Autgr S\subset GL_m\times GL_n$.
But the group on the right is contained in $Autgr S$ since it acts by
automorphisms of $S$ via linear changes of indeterminates.
It follows that $Autgr S\cong GL_m\times GL_n$.

Let $S=K(2k+1,n)=\prod_{j\geq -2}\g_j$ with its principal
grading.
By Proposition \ref{maximal}$(e)$, $Autgr S\subset 
\bar{\F}^\times\cdot(Sp_{2k}\times O_n)$. But the group on the right is 
contained in $Autgr S$ since it acts by automorphisms of $S$ via linear 
changes of indeterminates, preserving the supercontact differential form
$dx_{2k+1} + \sum_{i=1}^k(x_i dx_{k+i}-x_{k+i}dx_i)+\sum_{j=1}^n \xi_j d\xi_{n-j+1}$ up to multiplication
by a non-zero number.   
It follows that
$Autgr S\cong \bar{\F}^\times\cdot(Sp_{2k}\times O_n)$. 
Likewise, if $S=H(2k,n)=\prod_{j\geq -1}\g_j$ with 
the principal grading, the group
$\bar{\F}^\times\cdot(Sp_{2k}\times O_n)$ acts by automorphisms of $S$ via
linear changes of indeterminates, preserving the supersymplectic differential
form $\sum_{i=1}^k dp_i \wedge dq_i + \sum_{j=1}^n d\xi_j d\xi_{n-j+1}$ up to multiplication by
a non-zero number. Hence again,
$Autgr S\cong \bar{\F}^\times\cdot(Sp_{2k}\times O_n)$. 

Consider the Lie superalgebras  $S=HO(n,n)$,
$SHO(n,n)$ with $n>3$, $KO(n,n+1)$, or
$SKO(n,n+1;\beta)$ with $n>2$ or with $n=2$ and  $\beta\neq 0,1$,
 with  their principal
gradings.
By Proposition
\ref{maximal}$(f)$, $Autgr S\subset\bar{\F}^\times\cdot GL_n$. 
On the other hand, the group on the right is
contained in $Autgr S$ since it acts by automorphisms of $S$ via linear 
changes of indeterminates, preserving the odd
supersymplectic form $\sum_{i=1}^ndx_id\xi_i$ and 
the volume form $v$ attached to the usual divergence
up to multiplication
by a non-zero number, if $S=HO(n,n)$ or $S=SHO(n,n)$ with $n>3$,
and preserving the odd supercontact form
$d\xi_{n+1}+\sum_{i=1}^n(\xi_idx_i+x_id\xi_i)$
up to multiplication by an invertible function,
 and the volume form $v_{\beta}$ attached to 
the $\beta$-divergence
up to multiplication by a non-zero number, if $S=KO(n,n+1)$ or
$S=SKO(n,n+1;\beta)$ with $\beta\neq 0,1$ if $n=2$.   
Therefore $Autgr S\cong\bar{\F}^\times\cdot GL_n$.

Let $S=SHO^\sim(n,n)$, with $n>2$ even, 
and let $S_0$ be the principal subalgebra of $S$.
The group $Aut(GrS)$ 
consists of the automorphisms of $SHO'(n,n)$
preserving its principal grading. By the same
argument as for $SHO(n,n)$, $Aut(Gr S)\cong \bar{\F}^\times\cdot GL_n$.
The subgroup $Autgr S$ consists of the elements 
in $Aut(Gr S)$ which can be lifted to automorphisms
of $S$. Every element in $SL_n$ can be lifted to
an automorphism of $S$, since it preserves the form 
$Fv$ defining the Lie superalgebra $S$.
Besides, such automorphisms are inner and
act on $S$ via linear changes of variables.
On the contrary, the outer automorphisms of $SHO'(n,n)$ 
  do not preserve the form $Fv$ for any $t\in\bar{\F}$, hence they cannot be lifted to any
automorphism of $S$. It follows that $Autgr S\cong SL_n$. The argument for
$S=SKO^\sim(n,n+1)$ is similar.


Consider $S=S(1,2)=\prod_{j\geq -2}
\g_j$ with 
the grading of type $(2|1,1)$. 
By Proposition \ref{maximal}$(a)$, $Autgr S\subset\bar{\F}^\times\cdot SO_4$.
Notice that $\exp(ad(\g_0))
\cong \bar{\F}^\times\cdot SL_2$, acting by
automorphisms of $S$ via linear changes of indeterminates
which
preserve the standard volume form $v$ up to multiplication by a non-zero
number. 
Since the algebra of outer derivations of $S$ is isomorphic to $sl_2$,
$Autgr S\cong\bar{\F}^\times\cdot SO_4$.

Consider the Lie superalgebra $S=SHO(3,3)=\prod_{j\geq -j}\g_j$ with 
 the grading
of type $(2,2,2|1,1,1)$. 
By Proposition \ref{maximal}$(a)$,
$Autgr S\subset\bar{\F}^\times\cdot (SL_3\times SL_2)$.
Notice that $\exp(ad(\g_0)_{\bar{0}})\cong SL_3$,
acting by
automorphisms of $S$ via linear changes of indeterminates,
which preserve the odd
supersymplectic form  and 
the volume form $v$
up to multiplication
by a non-zero number.
Since the algebra of outer derivations of $S$ is isomorphic to
$gl_2$,  $Autgr S\cong \bar{\F}^\times\cdot (SL_3\times SL_2)$.

Consider the Lie superalgebra $S=SKO(2,3;1)=\prod_{j\geq -1}\g_j$,
with  its subprincipal grading.
In this case $\g_{-1}=\bar{\F}^{2|2}$ hence,
by Proposition \ref{Victor}$(a)$, $Autgr S\subset GL_2\times GL_2$.
As we recalled in Section \ref{gradings},
$(\g_0)_{\bar{0}}\cong sl_2+\bar{\F}$, where $sl_2$
acts trivially on $(\g_{-1})_{\bar{1}}$ and by the
standard action on $(\g_{-1})_{\bar{0}}$, and $\bar{\F}$ acts as the scalar $-2$ on $\g_{-2}$.
Besides, the algebra of outer derivations of $S$ is
isomorphic to $sl_2$, it acts trivially on $(\g_{-1})_{\bar{0}}$ and by the standard
action on $(\g_{-1})_{\bar{1}}$.
It follows that
 $Autgr S\cong \bar{\F}^\times\cdot (SL_2\times
SL_2)$.

Consider the Lie superalgebra $S=E(1,6)$ with 
its principal grading.
By Proposition \ref{maximal}$(e)$ with $k=0$ and $n=6$, 
$Autgr S\subset\bar{\F}^\times\cdot O_6$. Notice that
$\exp(ad(\g_0)_{\bar{0}})\cong \bar{\F}^\times \cdot SO_6$ 
and the group $O_6/SO_6$ is generated by the change of
variables $\xi_i\leftrightarrow \eta_i$ which is not an
automorphism of $E(1,6)$, since it exchanges the submodules
$\g_1^+$ and $\g_1^-$ of the $1$-st graded component of $K(1,6)$
in its principal grading (cf.\ \cite[\S 6]{CantaK}). Therefore $Autgr S\cong \bar{\F}^\times\cdot SO_6$.

Consider $S=E(3,6)=\prod_{j\geq -2}\g_j$ with its principal grading. 
By Proposition \ref{maximal}$(b)$, 
$Autgr S\subset \bar{\F}^\times\cdot (SL_2\times SL_3)$. Since
$\exp(ad(\g_0)_{\bar{0}})\cong \bar{\F}^\times\cdot (SL_2\times SL_3)$, the statement follows.
The same argument holds for $S=E(3,8)$ in its principal grading.

Consider $S=E(4,4)=\prod_{j\geq -1}\g_j$ with its principal grading.
By Proposition \ref{maximal}$(g)$, 
$Autgr S\subset \bar{\F}^\times\cdot SL_4$. Since
 $\exp(ad(\g_0)_{\bar{0}})\cong
GL_4$, equality holds. 

Finally, consider $S=E(5,10)=\prod_{j\geq -2}\g_j$ with its principal grading.
 By Proposition
\ref{maximal}$(c)$, $Autgr S\subset\bar{\F}^\times\cdot SL_5$.
Besides,
 $\exp(ad(\g_0))\cong SL_5$. Note that the Lie superalgebra $S$ has
an outer derivation acting on $S=\prod_{j\geq -2}\g_j$ as the
grading operator.
It follows that  $Autgr S\cong GL_5$. \hfill$\Box$

\begin{corollary}\label{outer} Let $S$ be a simple infinite-dimensional
linearly compact Lie superalgebra over $\bar{\F}$. Then $Aut S$ is the semidirect product
of the group of inner automorphisms of $S$  and the 
finite-dimensional
algebraic group $A$, described below:
\begin{itemize}
\item[$(a)$] if 
$S=SHO^\sim(n,n)$ or $SKO^\sim(n,n+1)$,
then $A=\{1\}$;
\item[$(b)$] if $S$ is a Lie algebra or $S=E(1,6)$, $E(3,6)$, $E(3,8)$,
$E(4,4)$,
$E(5,10)$, then $A\cong\bar{\F}^\times$;
\item[$(c)$] if $S=H(m,n)$ or $K(m,n)$ with $n>0$, 
then $A\cong\bar{\F}^\times\times\mathbb{Z}_2$;
\item[$(d)$] if $S=W(m,n)$ with $(m,n)\neq (1,1)$ and $n>0$,
 $S(m,n)$ with $m>1$ and $n>0$ or with $m=1$ and $n$ odd, $HO(n,n)$, 
$SHO(n,n)$ with $n>3$ even, $KO(n,n+1)$,
 $SKO(n,n+1;\beta)$ with $\beta\neq (n-2)/n,1$, or with $n$ even and
$\beta=(n-2)/n$, or with $n$ odd and $\beta=1$, 
then $A\cong \bar{\F}^{\times 2}$;
\item[$(e)$] if $S=S(1,n)$, with $n>2$ even, $SKO(n,n+1;(n-2)/n)$
with $n>2$ odd, $SKO(n,n+1;1)$ with $n>2$ even,
or $SHO(n,n)$ with
$n>3$ odd, 
then $A\cong U\rtimes\bar{\F}^{\times 2}$ where $U$ is a one-dimensional unipotent group;
\item[$(f)$] if $S=S(1,2)$,  then $A\cong \bar{\F}^\times\times SO_3$;
\item[$(g)$] if $S=SHO(3,3)$ or $SKO(2,3;1)$,  then $A\cong \bar{\F}^{\times}\cdot SL_2$.
\end{itemize}
\end{corollary}

\medskip

We shall now investigate the nature
of all continuous automorphisms of each simple
infinite-dimensional non-exceptional linearly compact Lie superalgebra $S$
over $\bar{\F}$.

\begin{lemma}\label{cv} Consider a subalgebra $L$
of $W(m,n)$ and let $D$ be an even
element of $L$ lying in the first member of   
a filtration of $W(m,n)$.
Then $D$ lies in the Lie algebra of the group of changes of variables
which map $L$ to itself. 
\end{lemma}
{\bf Proof.} Let $D =\sum_i P_i \frac{\partial}{\partial x_i}+\sum_j Q_j 
\frac{\partial}{\partial \xi_j}$. Then $\exp tD$, when 
applied to $x_i$ and $\xi_j$, gives  convergent series $S_i(t)$ and $R_j(t)$, 
respectively (in the linearly compact topology), hence the change of variables
$x_i \rightarrow S_i(t)$, $\xi_j\rightarrow R_j(t)$
 is a one-parameter group of automorphisms of
$W(m,n)$ which preserves $L$. 
\hfill$\Box$

\begin{theorem}\label{changesof var}
Let $S\subset W(m,n)$ be the defining embedding of a 
non-exceptional simple
infinite-dimensional linearly compact Lie superalgebra.
If $S\neq S(1,2)$, $SHO(3,3)$, $SKO(2,3;1)$, and
$S$ is defined by an action on a volume form $v$, an even
or odd supersymplectic form $\omega_s$, or an even or odd
supercontact form $\omega_c$, then all continuous automorphisms
of $S$ over $\bar{\F}$ are obtained by invertible changes of variables,
multiplying $v$ and $\omega_s$ by a constant
and $\omega_c$ by a function.
If $S=S(1,2)$, $SHO(3,3)$, or $SKO(2,3;1)$, then all these
changes of variables form a subgroup $H$ of $Aut S$ of
codimension one.
\end{theorem}
{\bf Proof.} Let $S=W(m,n)$ with $(m,n)\neq (1,1)$. Then $Der S=S$ hence,
by Lemma \ref{cv}, $Autf S$ consists of invertible changes of
variables. Besides, $Autgr S$ consists of linear changes of
variables, thus, by Proposition \ref{Victor}$(c)$, the statement
for $W(m,n)$
follows.

Let now $S=S(m,n)$ with $(m,n)\neq (1,2)$.
Then $Der S=CS'(m,n) \subset W(m,n)$. By Lemma \ref{cv}, $Autf S$ lies in the group of
changes of variables whose Lie algebra kills the volume form $v$
attached to the standard divergence. Hence
these changes of variables preserve the form $v$.
It is clear that all linear changes of variables multiply the
volume form $v$ 
by a number, hence all of them are automorphisms of $S(m,n)$.
Hence $Aut S$ is the group of changes of variables which preserve the volume 
form up to multiplication by a non-zero  number.

If $S=S(1,2)$, then, by the same argument as above, 
the inner automorphisms  of $S$ and its automorphism
$t^c$, where $c$ is the grading operator of $S$ with respect to its
grading of type $(2|1,1)$, are induced by
changes of variables which preserve the volume 
form up to multiplication by a non-zero  number. 
We recall that the algebra $\mathfrak{a}$ of outer derivations
of $S$ is isomorphic to $sl_2$, with standard generators
$e$, $f$, $h$, where $e=ad(\xi_1\xi_2\frac{\partial}
{\partial x})$ and $h=ad(\xi_1\frac{\partial}{\partial\xi_1}+
\xi_2\frac{\partial}{\partial\xi_2})$ (cf.\
\cite[Remark 2.12]{CantaK}).
 As for $S(m,n)$,
$t^h$, where $t\in\bar{\F}$, is obtained by
a linear change of variables preserving
the volume form up to multiplication by a non-zero number.
Besides, the element $\xi_1\xi_2\frac{\partial}
{\partial x}$ is contained in the first member of the
principal filtration of $S'(1,2)$, thus it is obtained
by a change of variables
preserving
the volume form up to multiplication by a non-zero number,
by Lemma \ref{cv}. On the other hand, the automorphism
$\exp(f)$ cannot be induced by any change of variables, since it
does not preserve the principal filtration of $S$.

%
The argument
for all other non-exceptional Lie superalgebras is similar.
\hfill$\Box$

\begin{remark}\em Let $S=S(1,2)$, $SHO(3,3)$, or $SKO(2,3;1)$.
In all these cases the algebra of outer derivations contains
$sl_2=\langle e,h,f\rangle$.
Denote by $U_{-}$ the
one-parameter
group of automorphisms $\exp(ad(tf))$,
where $f$ is explicitely described in
\cite[Remarks 2.12, 2.37, 4.15]{CantaK}.
Then $U_{-}$ is the  
``complementary" to $H$ subgroup in $Aut S$, namely,
for every $\varphi\in Aut S$, either $\varphi\in
U_{-}H$ or $\varphi\in U_{-}sH$ where $s$ is the reflection
$s=\exp(e)\exp(-f)\exp(e)$.
\end{remark}

\section{$\boldsymbol{\F}$-Forms}
Let ${\F}$ be a field of characteristic zero
and let $\bar{\F}$ be its algebraic closure.
\begin{definition} Let $L$ be a Lie superalgebra
over $\bar{\F}$. A Lie superalgebra $L^\F$ over $\F$ is called
an $\F$-form of $L$ if 
$L^\F\otimes_\F\bar{\F}\cong L$.
\end{definition}
Denote by $Gal$ the Galois group of $\F\subset \bar{\F}$.
Then $Gal$ acts on $Aut L$ as follows:
$$\alpha.\varphi:=\varphi^{\alpha}=\alpha\varphi\alpha^{-1}, ~~\alpha\in Gal, ~\varphi\in Aut L.$$
To any $\F$-form $L^\F$ of $L$, i.e., to any isomorphism
$\phi: {L^\F}\otimes_\F\bar{\F}\rightarrow L$, we can associate 
the map $\gamma_{\alpha}:Gal\rightarrow Aut L$,
 $\alpha\mapsto \phi^\alpha\phi^{-1}$.
The map $\gamma_{\alpha}$ satisfies the cocycle condition, i.e.,
$\gamma_{\alpha\beta}=\gamma_{\beta}^\alpha\gamma_{\alpha}$.
Two cocycles $\gamma$ and $\delta$ are equivalent if and only if
there exists an element $\psi\in Aut L$ such that
$\gamma_{\alpha}=(\psi^{-1})^\alpha\delta_{\alpha}\psi$. It follows
that equivalent cocycles correspond to isomorphic $\F$-forms.
\begin{proposition}\label{basic}
The map
$\phi\longmapsto \{\alpha\mapsto\phi^\alpha\phi^{-1}\}$
induces  a bijection between the set of isomorphism classes of
 $\F$-forms of $L$ and $H^1(Gal, Aut L)$.
\end{proposition}
{\bf Proof.} For a proof see \cite[\S 4]{R}.

\medskip

We recall the following standard result (cf.\ \cite[\S VII.2]{Serre1}):
\begin{proposition}\label{longexact}
If $K$ is a group and $A$, $B$, $C$ are groups with an action of $K$ by
automorphisms, related by
an exact sequence:
$$1\rightarrow A\rightarrow B\rightarrow C\rightarrow 1,$$
then there is a cohomology long exact sequence:
$$1\rightarrow H^0(K,A)\rightarrow H^0(K,B)\rightarrow H^0(K,C)
\rightarrow H^1(K,A)\rightarrow H^1(K,B)\rightarrow H^1(K,C),$$
where the first three maps are group homomorphisms, and the last three
are maps of pointed sets.
\end{proposition} 

\begin{proposition}\label{reduction} Let $S$ be a simple infinite-dimensional linearly
compact Lie superalgebra over $\bar{\F}$. Then
the map $j: Aut S \rightarrow Autgr S$ induces an embedding
$$j_{*}: H^1(Gal, Aut S)\longrightarrow H^1(Gal, Autgr S).$$
\end{proposition}
{\bf Proof.} The same arguments as in \cite[Proposition 4.2]{R}
show that $H^1(Gal,$ $Autf S)=0$. Then the statement follows 
from  exact sequence (\ref{exact}) in Section
\ref{general}  and Proposition \ref{longexact}. 
\hfill$\Box$

\bigskip

We recall the following well known results on 
Galois cohomology. All details can be found in 
\cite[\S ~X]{Serre1} and \cite[\S ~III Annexe]{Serre2}.
\begin{theorem}\label{results}
\begin{enumerate}
\item[$(a)$] $H^1(Gal, \bar{\F}^\times)=1$;
\item[$(b)$] $H^1(Gal, GL_n(\bar{\F}))=1$;
\item[$(c)$] $H^1(Gal, SL_n(\bar{\F}))=1$;
\item[$(d)$] $H^1(Gal, Sp_n(\bar{\F}))=1$;
\item[$(e)$] if $q$ is a quadratic form over $\F$, then
there exists a bijection between $H^1(Gal, O_n(q,\bar{\F}))$ and the set of
classes of $\F$-quadratic forms which are $\bar{\F}$-isomorphic to $q$;
\item[$(f)$] if $q$ is a quadratic form over $\F$, then
there exists a bijection between $H^1(Gal, SO_n(q,\bar{\F}))$ and the set of
classes of $\F$-quadratic forms ~$q'$ which are $\bar{\F}$-isomorphic to $q$ and
such that $\det(q')/\det(q)\in(\F^\times)^2$.
\end{enumerate}
\end{theorem}

\begin{lemma}\label{almostdirect}\em  Let $G$
   be an almost direct product over $\F$ of  
$\bar{\F}^\times$ and an algebraic group $G_1$, and let
$C=\bar{\F}^\times\cap G_1(\bar{\F})$ be a
cyclic group of order $k$. Then we have the following exact sequence:
\begin{equation}
1\rightarrow \F^\times/(\F^\times)^k\rightarrow H^1(Gal,G_1)\rightarrow
H^1(Gal,G)\rightarrow 1.
\label{!}
\end{equation}
In particular, if $H^1(Gal,G_1)=1$, then $H^1(Gal,G)=1$.
\end{lemma}
{\bf Proof.} 
We have the following exact sequence: 
$$1\rightarrow G_1\rightarrow G\stackrel{\pi}{\rightarrow}
\bar{\F}^\times\rightarrow 1,$$
where $\pi: G\rightarrow G/G_1\cong\bar{\F}^\times/C\rightarrow \bar{\F}^\times$
is the composition of the canonical map of $G$ to $G/G_1$
and the map $x\mapsto x^k$ from $\bar{\F}^\times/C$ to $\bar{\F}^\times$.
By Proposition
\ref{longexact}, we get the following exact sequence:
$$1\rightarrow \F^\times/C\rightarrow\F^\times\rightarrow H^1(Gal,G_1)\rightarrow
H^1(Gal,G)\rightarrow 1.$$
This implies exact sequence (\ref{!}).
\hfill$\Box$ 

\bigskip

We fix the $\F$-form $S^\F$ of each simple infinite-dimensional
linearly compact Lie superalgebra $S$ over $\bar{\F}$, defined by the
same conditions as in \cite{CantaK}, but over $\F$ (in the case
of $SKO(n,n+1;\beta)$ we need to assume that $\beta\in\F$).
This is called the {\em split} $\F$-form of $S$. In more invariant terms,
this $\F$-form is characterized by the condition that it contains
a split maximal torus $T$ (i.e.\ $T$ is $ad$-diagonalizable over $F$ and 
$T \otimes_\F\bar{\F}$
is a maximal torus of $S$).

\begin{theorem} Let $S$ be a simple infinite-dimensional linearly
compact Lie superalgebra  over $\bar{\F}$ not isomorphic to
$H(m,n)$, $K(m,n)$, $E(1,6)$, or $S(1,2)$.
Then any $\F$-form of $S$ is isomorphic to the
split $\F$-form. 
\end{theorem}
{\bf Proof.} It follows from Propositions \ref{basic}, \ref{reduction}
and the description of the group $Autgr S$ given in Theorem
\ref{AutgrS}, using Theorem \ref{results} and Lemma \ref{almostdirect}.
\hfill$\Box$

\begin{remark}\label{HeK}\em Let $S=H(2k,n)$ or
$S=K(2k+1,n)$. Then, according to Table 1 and 
Lemma \ref{almostdirect}, we have the exact sequence
$$1\rightarrow \F^\times/(\F^\times)^2\rightarrow H^1(Gal, Sp_{2k}\times O_n)\rightarrow
H^1(Gal, Autgr S)\rightarrow 1.$$
Here and further, $Sp_{2k}=Sp_{2k}(\bar{\F})$ and $O_n\subset
GL_n(\bar{\F})$ is the orthogonal group over $\bar{\F}$ which
leaves invariant the quadratic form $\sum_{i=1}^nx_ix_{n-i+1}$.
Since $H^1(Gal,$ $G_1\times G_2)\cong H^1(Gal, G_1)\times H^1(Gal, G_2)$,
by Theorem \ref{results}$(d)$, $H^1(Gal, Sp_{2k}\times O_n)\cong H^1(Gal,$ $O_n)$, hence we have 
the exact sequence
$$1\rightarrow \F^\times/(\F^\times)^2\rightarrow H^1(Gal, O_n)\rightarrow
H^1(Gal, Autgr S)\rightarrow 1.$$
\end{remark}

\medskip

Given a non-degenerate quadratic form $q$ over $\F$ in $n$ indeterminates, 
associated to a symmetric  matrix $c=(c_{ij})$,
introduce the 
following supersymplectic and supercontact differential forms $\sigma_q$
and $\Sigma_q$:
$$\sigma_q = \sum_{i=1}^k dp_i\wedge dq_i+ \sum_{i,j=1}^n c_{ij} d\xi_i d\xi_j,$$
$$\Sigma_q= dt+\sum_{i=1}^k (p_idq_i-q_idp_i)+\sum_{i,j=1}^n c_{ij} \xi_i d\xi_j.$$

\begin{theorem}\label{H(m,n)} 
$(a)$ Any $\F$-form of the Lie superalgebra
$S=H(2k,n)$ is isomorphic to one of  the Lie superalgebras
$H_q (2k,n):=\{X\in W(2k,n)^{\F}~|~X\sigma_q=0\}$. 

\noindent
$(b)$ Any $\F$-form of the Lie superalgebra
$S=K(2k+1,n)$ is isomorphic to one of the Lie superalgebras
$K_q (2k+1,n):=\{X\in W(2k+1,n)^{\F}~|~X\Sigma_q=f\Sigma_q\}$. 

Two such $\F$-forms
$S_q$ and $S_{q'}$ of $S$ are isomorphic if and only if
$q$ and $q'$ are equivalent non-degenerate quadratic forms over $\F$, up to multiplication by a
non-zero scalar in $\F$.
\end{theorem}
{\bf Proof.} It is easy to see that every non-degenerate
quadratic form $q$ over $\F$, with matrix $c=(c_{ij})$,
gives rise to the $\F$-forms $H_q(2k,n)$ and
$K_q(2k+1,n)$ of the Lie superalgebras
$S=H(2k,n)$ and $S=K(2k+1,n)$, respectively, attached to the corresponding
cocycles.
By construction, equivalent quadratic forms give rise to isomorphic
$\F$-forms of $S$. Besides, if $\lambda\in\F^\times$ and
$q'$ is the quadratic form associated to the matrix $\lambda c$, then
$S_q\cong S_{q'}$, 
and the isomorphism is given by the
following change of variables:
$$p_i\mapsto \lambda^{-1}p_i, ~~q_i\mapsto q_i, ~~\xi_i\mapsto \xi_i, 
~~~\mbox{if}~ S=H(2k,n)$$
$$t\mapsto\lambda^{-1}t, ~~p_i\mapsto \lambda^{-1}p_i, ~~q_i\mapsto q_i, ~~\xi_i\mapsto \xi_i, 
~~~\mbox{if}~ S=K(2k+1,n).$$
The $\F$-forms $S_q$ exhaust all $\F$-forms of the Lie superalgebra $S$,
due to Proposition \ref{basic}, Theorem \ref{AutgrS},
Remark \ref{HeK} and Theorem \ref{results}$(e)$.
\hfill$\Box$

\begin{example}\label{realE(1,6)}\em Consider the
$\F$-form $K_q(1,6)$ of $K(1,6)$ corresponding to the supercontact
form 
$\Sigma_q=dt+\sum_{i=1}^6c_{ij}\xi_id\xi_j$.
Then the principal grading of $K(1,6)$ induces an irreducible grading
on $K_q(1,6)$: $K_q(1,6)=\prod_{j\geq -2}\mathfrak{g}_j$,
where $\g_0=\h\oplus\F$, $\h$ is an $\F$-form of $so_6(\bar{\F})$, 
$\mathfrak{g}_{-1}\cong\F^6$, 
and  
$\g_1=\mathfrak{g}_{-1}^*\oplus \Lambda^3(\F^6)$. 

Let $d$ be the discriminant of the quadratic form $q$.
If $-d\in(\F^\times)^2$, then
the $\g_0$-module $\Lambda^3(\F^6)$ is not irreducible,
and  decomposes over $\F$ into
the direct sum of two $\g_0$-submodules
$\mathfrak{g}_1^+$ and $\mathfrak{g}_1^-$,
which are the eigenspaces of the Hodge operator ${}^*$, see
Example \ref{CK6} below
(they are obtained from one another by an
automorphism of $\g_0$).
 It follows that we can define an $\F$-form $E_q(1,6)$ of
the Lie superalgebra $E(1,6)$ by repeating the same construction
as the one described in Section \ref{gradings}, namely,
$E_q(1,6)$ will be the graded subalgebra of $K_q(1,6)$ generated by
$\mathfrak{g}_{-1}+\mathfrak{g}_0+(\mathfrak{g}_{-1}^*+\mathfrak{g}_1^+)$. 
\end{example}

\begin{theorem}\label{FformsofE(1,6)} Any $\F$-form of the Lie superalgebra
 $S=E(1,6)$ is isomorphic to one of the Lie superalgebras
$E_q(1,6)$ constructed in Example \ref{realE(1,6)}, where $q$ is
a non-degenerate quadratic form over $\F$ in six indeterminates,
with discriminant 
$d\in-(\F^\times)^2$.

Two such $\F$-forms
$E_q(1,6)$ and $E_{q'}(1,6)$ of $E(1,6)$ are isomorphic if and only if
the quadratic forms $q$ and $q'$
are equivalent, up to multiplication by 
a non-zero scalar in $\F$.
\end{theorem}
{\bf Proof.} By Lemma \ref{almostdirect} and
Theorem \ref{AutgrS}$(h)$, we have the exact sequence
$$1\rightarrow \F^\times/(\F^\times)^2\rightarrow H^1(Gal, SO_{6})\rightarrow
H^1(Gal, Autgr S)\rightarrow 1.$$
The statement follows,
due to Proposition \ref{basic},
Theorem \ref{results}$(f)$ and the proof of Theorem \ref{H(m,n)}.
\hfill$\Box$

\begin{remark}\label{S(1,2)inK(1,4)}\em
Consider the Lie superalgebra $K(1,4)=\prod_{j\geq -2}\g_j$
over $\bar{\F}$ with  respect to its principal grading. Then:
$\g_0\cong cso_4=so_4+\bar{\F}t$,
$\g_{-1}\cong\bar{\F}^4$ and
$\g_{-2}=[\g_{-1}, \g_{-1}]\cong\bar{\F}$.  Besides,
$\g_1=V_1\oplus V_{-1}$, where
for every $\lambda=\pm 1$, $V_{\lambda}$ is isomorphic to the
standard $so_4$-module, $[V_{\lambda}, V_{\lambda}]$ is
isomorphic to the trivial $so_4$-module $\bar{\F}$, and 
$\g_{-2}+\g_{-1}+gl_2+V_{\lambda}+[V_{\lambda},V_{\lambda}]\cong sl(2,2)/\bar{\F}$.
Finally,
$\g_2=\bar{\F}\oplus
so_4\oplus\bar{\F}$, where $so_4$ and $\bar{\F}$ denote the
adjoint and the trivial $so_4$-module, respectively.
Here $t$ acts as the grading operator.

The Lie superalgebra $S(1,2)$ over $\bar{\F}$ 
is the subalgebra of
$K(1,4)$ generated by $\h_{-1}\oplus \h_0\oplus \h_1\oplus\h_2$
where 
$\h_{-1}=\g_{-1}$, 
$\h_0=sl_2+\bar{\F}t$,
$\h_1=V_1$ and
$\h_2=\bar{\F}+sl_2$, where $sl_2$ denotes the adjoint $sl_2$-module
(see also \cite[Remark 2.33]{CantaK}).
\end{remark}
 
\begin{example}\label{S(1,2)}\em
Consider a Lie superalgebra $\g_{-}$ over $\F$ with
consistent $\Z$-grading $\g_{-}=\g_{-2}+\g_{-1}+\g_0$,
where $\g_0=\h\oplus\F$, $\h$ is an $\F$-form of $sl_2$,
$\g_{-2}=\F z$, and $\g_{-1}$ is a four-dimensional
$[\g_0, \g_0]$-module such that $\g_{-1}\otimes_\F\bar{\F}$ is
the direct sum of two copies of the standard $sl_2$-module,
and where the bracket in $\g_{-1}$ is defined as follows:
$$[a,b]=q(a,b)z,$$
where $q$ is a non-degenerate bilinear form on $\g_{-1}$
over $\F$, which is symmetric, i.e. $q(a,b)=q(b,a)$. 
Such a superalgebra $\g_{-}$ exists if and only if the
discriminant $d$ of the quadratic form $q$ lies in
$(\F^\times)^2$.
Let $S_q(1,2)$ denote the full prolongation of $\g_{-}$ over $\F$
(see \cite[\S 1.6]{CK}).
Then $S_q(1,2)$ is an $\F$-form of $S(1,2)$.
\end{example}

\begin{theorem}\label{Fformsof S(1,2)} Any $\F$-form of the Lie superalgebra
 $S=S(1,2)$ is isomorphic to one of  the Lie superalgebras
$S_q(1,2)$ constructed in Example \ref{S(1,2)}, where $q$
is a non-degenerate quadratic form over $\F$ in four indeterminates, with discriminant
$d\in(\F^\times)^2$.

Two such $\F$-forms
$S_q(1,2)$ and $S_{q'}(1,2)$ of $S(1,2)$ are isomorphic if and only if
the quadratic forms $q$ and $q'$ are equivalent, 
up to multiplication by a non-zero scalar in $\F$.
\end{theorem}
{\bf Proof.} By Lemma \ref{almostdirect} and
Theorem \ref{AutgrS}, we have the exact sequence
$$1\rightarrow \F^\times/(\F^\times)^2\rightarrow H^1(Gal, SO_{4})\rightarrow
H^1(Gal, Autgr S)\rightarrow 1.$$
The statement follows,
due to Proposition \ref{basic},
Theorem \ref{results}$(f)$ and the proof of Theorem \ref{H(m,n)}.
\hfill$\Box$

\bigskip

We summarize the results of this section in the following theorem:
\begin{theorem}\label{summary}
Let $S$ be a simple infinite-dimensional linearly
compact Lie superalgebra  over $\bar{\F}$. 
If $S$ is not isomorphic to $H(m,n)$, $K(m,n)$, $E(1,6)$, or $S(1,2)$, then
the split $\F$-form $S^\F$ is, up to isomorphism, the unique
$\F$-form of $S$. In the remaining four cases, all
$\F$-forms of $S$ are, up to isomorphism, as follows:
\begin{itemize}
\item[$(a)$] the Lie superalgebras
$H_q(m,n):=\{X\in W(m,n)^{\F}~|~X\sigma_q=0\}$ where $\sigma_q$
is a supersymplectic differential form over $\F$, if $S=H(m,n)$;
\item[$(b)$] the Lie superalgebras
$K_q(m,n):=\{X\in W(m,n)^{\F}~|~X\Sigma_q=f\Sigma_q\}$
where $\Sigma_q$ is a supercontact differential form over $\F$, if $S=K(m,n)$;
\item[$(c)$] the Lie superalgebras $E_q(1,6)$ constructed in
Example \ref{realE(1,6)}, where $q$ is a non-degenerate quadratic form
over $\F$ in six indeterminates
with discriminant $d\in-(\F^\times)^2$, if $S=E(1,6)$;
\item[$(d)$] the Lie superalgebras $S_q(1,2)$ constructed in
Example \ref{S(1,2)}, where $q$ is a non-degenerate quadratic form
over $\F$ in four indeterminates
with discriminant $d\in(\F^\times)^2$, if $S=S(1,2)$.  
\end{itemize}
The isomorphisms between these $\F$-forms are described in Theorems
\ref{H(m,n)}, \ref{FformsofE(1,6)},
\ref{Fformsof S(1,2)}.
\end{theorem}

\begin{remark}\label{real}\em It follows immediately from
Theorem \ref{summary} that a
simple infinite-dimensional linearly
compact Lie superalgebra $S$  over $\C$ has, up to isomorphism,
one real form 
if $S$ is not isomorphic to $H(m,n)$, $K(m,n)$, $E(1,6)$, or $S(1,2)$, two real forms
if $S$ is isomorphic to $E(1,6)$ or $S(1,2)$, and $[n/2]+1$ real forms if
$S$ is isomorphic to $H(m,n)$ or $K(m,n)$.
\end{remark}

\section{Finite Simple Lie Conformal Superalgebras}
In this section we use the theory of Lie conformal superalgebras
in order to give an explicit construction of all non-split forms of all simple
infinite-dimensional linearly compact Lie superalgebras. 
In conclusion of the section, we give the related classification of all
$\F$-forms of all simple finite Lie conformal superalgebras.

We briefly recall the definition of a Lie conformal superalgebra
and of its annihilation algebra. For notation, definitions and results
on Lie conformal superalgebras we refer to \cite{DK}, \cite{FK}
and \cite{V}.

A Lie conformal superalgebra $R$ over $\F$
is a left $\Z/2\Z$-graded $\F[\partial]$-module endowed with an $\F$-linear
map, called the $\lambda$-bracket,
$$R\otimes R\rightarrow \F[\lambda]\otimes R, ~a\otimes b\mapsto[a_\lambda b],$$
satisfying the axioms of sesquilinearity, skew-commutativity, and the
Jacobi identity.
One writes $[a_\lambda b]=\sum_{n\in\Z_+}\frac{\lambda^n}{n!}(a_{(n)}b)$;
the coefficient $(a_{(n)}b)$ is called the $n$-th product of $a$ and
$b$.
A Lie conformal superalgebra $R$ is called {\em finite} if it is finitely
generated as an $\F[\partial]$-module.

Given a finite Lie conformal superalgebra $R$, we can associate to it a 
linearly compact Lie superalgebra $L(R)$ as follows. Consider the
Lie conformal superalgebra $R[[t]]$, where $t$ is an even indeterminate, 
the $\partial$-action is defined by 
$\partial+\partial_t$, 
and the $n$-th products are defined by:
$$a(t)_{(n)}b(t):=\sum_{j\geq 0}(\partial_t^j a(t))_{(n+j)}b(t)/j!,$$
where $a(t)$,$b(t)\in R[[t]]$, 
and the $n$-th products on the right are extended from $R$ to $R[[t]]$
by bilinearity. 
Then $(\partial+\partial_t)R[[t]]$ is a two-sided ideal of $R[[t]]$
 with respect 
to $0$-th product, and this product induces a Lie superalgebra bracket on
$L(R):=R[[t]]/(\partial +\partial_t)R[[t]]$.
The linearly compact Lie superalgebra $L(R)$ is called the 
{\em annihilation algebra} of $R$.

For the classification of finite simple Lie conformal superalgebras over
an algebraically closed field $\bar{F}$ of characteristic zero  
we refer to \cite{FK}. The list consists of four series ($N\in\Z_+$):
$W_N$, $S_{N+2,a}$, $\tilde S_{N+2}$, $K_N$ ($N \neq 4$), $K'_4$,
the exceptional Lie conformal superalgebra $CK_6$ of rank 32, and
$Cur \mathfrak{s}$, where $\mathfrak{s}$ is a simple finite-dimensional
Lie superalgebra.

\begin{example}\label{K_N}\em
Let $V$ be an $N$-dimensional vector space over $\F$ with a
non-degenerate symmetric bilinear
form $q$. The Lie conformal superalgebra $K_{N,q}$ 
associated to $V$ is $\F[\partial]\Lambda(V)$ with $\lambda$-bracket:
\begin{equation}
[A_{\lambda} B]=(\frac{r}{2} -1)\partial(AB) +
(-1)^r \frac{1}{2} \sum_{j=1}^N (i_{a_j}A)(i_{b_j}B)
+\lambda (\frac{r+s}{2} -2)  AB,
\label{lambdabracket}
\end{equation}
where $A,B \in \Lambda (V)$, $r=\deg(A)$, $s=\deg(B)$, $a_i,b_i \in V$,
$q(a_i,b_j)=\delta_{i,j}$, and $i_a$, for $a \in V$, denotes the contraction
with $a$, i.e., $i_a$ is the odd derivation of  $\Lambda(V)$
defined by: $i_a(b)=q(a,b)$ for $b \in V$ (cf \cite[Example 3.8]{FK}). 
The annihilation algebra of $K_{N,q}$ is isomorphic to the Lie superalgebra
$K_q(1,N)$ defined in Theorem \ref{H(m,n)}$(b)$.
The Lie conformal superalgebra $K_{N,q}$ is an $\F$-form of the finite 
simple Lie conformal
superalgebra $K_N$. 
\end{example}

\begin{example}\label{CK6}\em Let $V$
be an $N$-dimensional vector space over $\F$ with a non-degenerate symmetric
bilinear form $q$, and let $K_{N,q}$ be the Lie
conformal superalgebra  over $\F$ constructed in
Example \ref{K_N}. Choose a basis $\xi_1, \dots, \xi_N$ of $V$, and let
 ${}^*$ denote the Hodge star
operator on $V$ associated to the form $q$, i.e., 
$$(\xi_{j_1} \wedge \xi_{j_2} \wedge\dots \wedge\xi_{j_k})^*=
i_{\xi_{j_1}}i_{\xi_{j_2}}\dots i_{\xi_{j_k}} 
(\xi_1\wedge\dots \wedge \xi_N).$$
It is easy to check that, for every $a\in\Lambda(V)$,
 $(a^*)^*=(-1)^{N(N-1)/2}\det(q)a$.

Let $N=6$ and choose $\alpha\in \F$ such that
$\alpha^2=-1/det(q)$. 
Consider the following elements in $\Lambda(V)$:
$$-1+\alpha\partial^3 1^*,  ~\xi_i\xi_j+\alpha\partial(\xi_i\xi_j)^*,
~\xi_{i}-\alpha\partial^2\xi_{i}^*, ~\xi_i\xi_j\xi_k+\alpha(\xi_i\xi_j\xi_k)^*.$$
It is easy to check that 
the $\F[\partial]$-span of these elements
is closed under $\lambda$-bracket (\ref{lambdabracket}),
hence they form an $\F$-form $CK_{6,q}$ of the Lie
conformal subalgebra $CK_6$ of $K_6$ (cf \cite[Theorem 3.1]{ChengK}).

Likewise, if $N=4$ and $\beta^2=1/det(q)$,
the $\F[\partial]$-span of the elements:
$$-1-\beta\partial^2 1^*, ~\xi_i\xi_j-\beta(\xi_i\xi_j)^*,
~\xi_{i}+\beta\partial\xi_{i}^*$$
is closed under $\lambda$-bracket (\ref{lambdabracket}).
It follows that these elements form a 
 subalgebra $S_{2,q}$ of  $K_{4,q}$, which is an $\F$-form of the
Lie conformal superalgebra $S_{2,0}$ (cf. \cite[Remark p.\ 225]{ChengK}).
\end{example}

\begin{remark}\label{explicit}\em
The annihilation algebras of the Lie conformal superalgebras 
$CK_{6,q}$ and $S_{2,q}$, constructed in Examples \ref{K_N} and
  \ref{CK6}, are
the Lie superalgebras  $E_q(1,6)$ and $S_q(1,2)$, constructed
in Examples \ref{realE(1,6)} and \ref{S(1,2)}, respectively.
Due to Theorems \ref{FformsofE(1,6)} and
\ref{Fformsof S(1,2)},
Example \ref{CK6} provides an explicit
construction of all $\F$-forms of the Lie superalgebras
$E(1,6)$ and $S(1,2)$.
\end{remark}

We conclude by classifying all $\F$-forms of all simple finite
Lie conformal superalgebras over 
$\bar{F}$. The following theorem can be derived
from \cite[Remark 3.1]{FKR}.

\begin{theorem}\label{conformalsummary}
Let $R$ be a simple finite Lie conformal superalgebra over
$\bar{\F}$. If $R$ is not isomorphic to $S_{2,0}$, $K_N$, $K'_4$, $CK_6$ or
$Cur \mathfrak{s}$, then there exists, up to isomorphism, a unique
$\F$-form of $R$ (in the case $R=S_{N,a}$, we have to assume that $a\in\F$
for such a form to exist). In the remaining  cases, all
$\F$-forms of $R$ are as follows:
\begin{itemize}
\item the Lie conformal superalgebras $K_{N,q}$ if $R$ is isomorphic
to $K_N$;
\item the Lie conformal superalgebras $CK_{6,q}$
 if $R$ is isomorphic to $CK_6$;
\item the Lie conformal superalgebras
$S_{2,q}$  if $R$ is isomorphic to $S_{2,0}$;
\item the derived algebras of the Lie conformal superalgebras
$K_{4,q}$ if $R$ is isomorphic to $K'_4$;
\item the Lie conformal superalgebras
$Cur {\mathfrak s}^\F$, where $\mathfrak{s}^\F$ is an $\F$-form of the Lie
superalgebra $\mathfrak{s}$, if $R$ is isomorphic to $Cur \mathfrak{s}$.
\end{itemize}
\end{theorem}

\medskip

\section*{Acknowledgments}
We would like to acknowledge the help of A.\ Maffei and E.\ B.\ Vinberg.

$$$$ 


\medskip

\begin{thebibliography}{99} 

%
%
%
%
%
\bibitem{CCK} 
{\sc N.\ Cantarini, S.-J.\ Cheng, V.\ G.\ Kac} 
\newblock{\em Errata: Structure of some 
$\mathbb Z$-graded Lie superalgebras of vector fields,} 
\newblock Transf. Groups {\bf 9} (2004), 399-400. 

\bibitem{CantaK} 
{\sc N.\ Cantarini, V.\ G.\ Kac} 
\newblock{\em Infinite dimensional primitive linearly compact Lie
superalgebras,} 
\newblock math.\ QA/0511424. 


\bibitem{ChengK} 
{\sc S.-J.\ Cheng, V.\ G.\ Kac} 
\newblock{\em 
A new $N=6$ superconformal algebra,} 
\newblock Comm.\ Math.\ Phys.\ {\bf 186} (1997), 219--231. 

\bibitem{CK} 
{\sc S.-J.\ Cheng, V.\ G.\ Kac} 
\newblock{\em Structure of some 
$\mathbb Z$-graded Lie superalgebras of vector fields,} 
\newblock Transf.\ 
Groups {\bf 4}  (1999), 219--272. 

\bibitem {DK}
{\sc A.\ D'Andrea, V.\ G.\ Kac} 
\newblock{\em Structure theory of finite conformal algebras,} 
\newblock Selecta Math.\ (N.S.) {\bf 4} (1998), 377--418. 

\bibitem{FK}
{\sc D.\ Fattori, V.\ G.\ Kac}
\newblock {\em Classification of finite simple Lie conformal superalgebras,}
\newblock J.\ Algebra {\bf 258} (2002), 23--59.

\bibitem{FKR}
{\sc D.\ Fattori, V. G. Kac, A. Retakh}
\newblock {\em Structure theory of finite Lie conformal superalgebras,}
\newblock In ``Lie theory and its applications to physics'', eds H.-D.\ 
Doebner
and V.\ V.\ Dobrev, World Sci., (2004) 27--63.

%
\bibitem{G2}
{\sc V.\ W.\ Guillemin}
\newblock{\em 
Infinite-dimensional primitive Lie algebras,}
\newblock J.\ Diff.\ Geom.\ {\bf 4} (1970), 257--282.
%
\bibitem{K2}
{\sc V.\ G.\ Kac}
\newblock {\em Lie superalgebras,}
\newblock Adv.\ Math.\ {\bf 26} (1977), 8--96.

%
\bibitem{K} 
{\sc V.\ G.\ Kac} 
\newblock{\em Classification of infinite-dimensional simple linearly 
compact Lie superalgebras,} 
\newblock Adv
.\ Math.\ {\bf 139} (1998), 1--55. 

\bibitem{V} 
{\sc V.\ G.\ Kac} 
\newblock{\em Vertex Algebras for Beginners,} 
\newblock Second Edition, University Lecture Series {\bf 10}, 
 AMS, 1998. 

\bibitem{Gafa} 
{\sc V.\ G.\ Kac} 
\newblock{\em Classification of infinite-dimensional simple groups of
supersymmetries and quantum field theory,}
\newblock GAFA, Geom.\ Funct.\ Anal.\ Special Volume GAFA2000 (2000), 162--183.
%

\bibitem{R}
{\sc A.\ N.\ Rudakov}
\newblock{\em Groups of automorphisms  of infinite-dimensional  simple 
Lie algebras,}
\newblock Math.\ USSR-Izvestija Vol.\ {\bf 3(4)} (1969), 707--722.

\bibitem{Serre1}
{\sc J.P.\ Serre}
\newblock{\em Corps locaux,}
\newblock Hermann, Paris, 1962.

\bibitem{Serre2}
{\sc J.P.\ Serre}
\newblock{\em Cohomologie Galoisienne. Cinqui\`eme \'edition,}
\newblock Lecture Notes in Mathematics {\bf 5}, Springer-Verlag, 
Berlin-Heidelberg, 1994.


\bibitem{S} {\sc I.\ Shchepochkina} 
\newblock{\em The five exceptional simple 
Lie superalgebras of vector fields and their fourteen regradings,} 
\newblock Repr.\ Theory 3 (1999), 373--415.
%

\end{thebibliography}
\end{document}